\documentclass[12pt]{article}
\usepackage{color,amssymb,amsthm,amsmath,graphicx,hyperref,float,bm,array,tikz}
\usepackage[headings]{fullpage}

\usepackage{fancyhdr}
\usepackage{xcolor}
\definecolor{dark-blue}{rgb}{0,0,0.45}
\hypersetup{colorlinks, linkcolor={dark-blue}, citecolor={blue}, urlcolor={blue}}

\numberwithin{equation}{section}

\def\fl#1{\left\lfloor#1\right\rfloor}
\def\Z{\mathbb{Z}}
\def\C{\mathbb{C}}

\newtheorem{theorem}{Theorem}[section]
\newtheorem{lemma}[theorem]{Lemma}

\newtheorem{corollary}[theorem]{Corollary}
\theoremstyle{definition}
\newtheorem{example}[theorem]{Example}
\newtheorem{definition}[theorem]{Definition}
\newtheorem{construction}[theorem]{Construction}

\title{A Recursive Construction of Perfect $K^n\times K^n$ Arrays\\ over $K$ Roots of Unity}
\author{Sam Blake}
\date{7 October 2026}

\begin{document}

\maketitle

\begin{abstract}
We give a uniform, explicit and elementary recursive construction of two-dimensional arrays with
perfect (ideal) periodic autocorrelation over a \emph{fixed} alphabet of $K$ roots of unity and
\emph{unbounded} size: for every integer $K\ge2$ and every $n\ge1$ a perfect $K^n\times K^n$ array
whose entries are $K$-th roots of unity, given by a nested sum of floor terms. The construction
works with \emph{twisted} (quasi-perfect) arrays indexed by a type $v\in\Z_K^2$, of which ordinary
perfect arrays are the type $v=0$ members. For $K=2$ the types are Jedwab's generalised perfect
binary arrays and the recursion specialises to the doubling constructions of Jedwab--Mitchell and
Wild; for prime $K$ the arrays are the relative difference sets whose existence in abelian
$p$-groups was settled by Ma and Schmidt; and the bare existence of perfect $K^n\times K^n$
arrays over $K$ roots of unity follows, for every $K$, from recent constructions of group-invariant
Butson Hadamard matrices. The contribution of this note is therefore not existence but the
calculus: types with arbitrary entries in $\Z_K^2$; three spectral lemmas (a twisted periodic
extension of a type-$v$ array has a flat spectrum supported on one coset of $K\Z^2$; interleaving
arrays whose spectra are flat on disjoint cosets gives a perfect array; and a shear of the index plane with a floor-function
carry permutes types through the action of $SL_2(\Z_K)$ on $\Z_K^2$, whose orbits are the classes
$\gcd(v_1,v_2,K)=g$); a one-directional lift that produces every type at side $KM$ from every type
at side $M$ for arbitrary $K$; and the resulting closed index function, which is the recursive form
of the floor-function arrays of \cite{Blake2017}.
\end{abstract}

\section{Introduction}

Let $\omega=e^{2\pi\sqrt{-1}/K}$. A two-dimensional array $S=[S_{x_1x_2}]$ of size $L_1\times L_2$ with
entries in $\mu_K=\{1,\omega,\dots,\omega^{K-1}\}$ is \emph{perfect} if all its off-peak periodic
autocorrelations vanish. The \emph{phase} (alphabet size) is $K$. A basic question, raised in
\cite{Blake2017} and in a different language by Mow \cite{Mow1993,Mow1996}, is how large a
perfect array can be at a fixed phase. A single floor term $\omega^{\lfloor x_1x_2/q\rfloor}$ reaches
side $2K^2$ for even $K$ (Construction~V of \cite{Blake2017}) and side $K^2$ for odd $K$ (the
Frank--Heimiller array \cite{Heimiller1961,Frank1962}), and one is led to conjecture that no single floored product
$\omega^{\lfloor p(x_1,x_2)/D\rfloor}$ exceeds side $2K^2$. The present note shows
that if one allows \emph{nested} floor terms the side is unbounded at every phase: for every integer
$K\ge2$ and every $n\ge1$ we give an explicit perfect $K^n\times K^n$ array over the $K$-th roots of
unity (Theorem~\ref{mainthm} in Section~\ref{towersec}). The smallest instances beyond the
single-floor family are a perfect $16\times16$ binary array, a perfect $27\times27$ array over the
cube roots of unity, and a perfect $125\times125$ array over the fifth roots of unity; for odd prime
$K$ the side $K^3$ already exceeds the single-floor ceiling $K^2$ by a factor $K$.

We stress at the outset that this is not new as an \emph{existence} statement.
For $K=2$ it is the square case of the Turyn--Davis--Kraemer theorem on Menon difference sets in
abelian $2$-groups \cite{Turyn1965,Davis1991,Kraemer1993} (the square case itself is
\cite[Thm.~4.2]{Davis1991}), realised explicitly by the recursive
array constructions of Jedwab and Mitchell and of Wild \cite{JedwabMitchell1988,Wild1988,
JedwabMitchell1990,Jedwab1992} and by Dillon's closed formula \cite{Dillon1990}. For $K=p$ prime it is
contained in the theorem of Ma and Schmidt \cite{MaSchmidt1995} that a
$(p^{2c},p,p^{2c},p^{2c-1})$-relative difference set exists in an abelian group of order $p^{2c+1}$
if and only if the exponent is at most $p^{c+1}$, and more recently in Do Duc's local-ring
construction of group-invariant Butson Hadamard matrices \cite{DoDuc2021}. For prime powers
$K=p^e$ a $p$-ary perfect array is a fortiori an array over $K$-th roots, and for composite $K$ the
prime-power cases combine by the Chinese remainder theorem (L\"uke's product, \cite{Luke1988}).
Section~\ref{litsec} gives the details. What we believe to be new is the uniform explicit formula and
the machinery in which it is expressed: the full set of $K^2$ types, the $SL_2(\Z_K)$ action on
them, and a lift theorem valid for every $K$. The closed formula is a $K$-ary analogue of Dillon's
explicit rank-two formula for Menon difference sets in $\Z_{2^{a+1}}^2$ \cite{Dillon1990}, quoted in
\cite[Thm.~4.1]{DavisJedwab1996}; it does not bear on the open problem posed there, which asks for
explicit constructions in abelian $2$-groups of higher rank.

The idea behind the recursion is as follows. A perfect array of side $KN$ decomposes into $K^2$
polyphase components of side $N$, and it is perfect precisely when their spectra are flat on
disjoint sets that cover the frequency plane. One would like to build such an array from
\emph{perfect} components, but an elementary counting argument (Section~\ref{lemmasec}) shows that
two perfect arrays can never have disjoint spectra, so the naive recursion ``perfect $\to$ perfect''
through components is impossible; this is precisely why the Frank--Heimiller array $\omega^{\lfloor x_1x_2/K\rfloor}$ and its
relatives, whose components are constants times a Frank--Heimiller core, stop at side $K^2$.
The components must instead be \emph{twisted} tilings of smaller arrays which are perfect only after
multiplication by a fractional linear phase.
These twisted arrays, which we call arrays of type $v\in\Z_K^2$, are closed under the recursion
provided one tracks all $K^2$ types simultaneously. For $K=2$ they are exactly the generalised perfect
binary arrays of Jedwab \cite{Jedwab1992}, and for general $K$ with $v\in\{0,1\}^2$ they are the
generalised perfect $K$-phase arrays of Armario, Egan and Flannery \cite{AEF2024}; the point of
this note is that with all of $\Z_K^2$ as the type space the same machinery, made explicit, works
verbatim for every $K$.

\section{Definitions and the spectral criterion}\label{defsec}

Throughout, arrays are indexed by $\Z_{L_1}\times\Z_{L_2}$, i.e.\ by $0\le x_i<L_i$ with indices read
modulo $L_i$. For a complex array $A$ on $\Z_{L_1}\times\Z_{L_2}$ write
$$\widehat A(\xi_1,\xi_2)=\sum_{x_1,x_2}A(x_1,x_2)\,e^{-2\pi\sqrt{-1}\,(\xi_1x_1/L_1+\xi_2x_2/L_2)},
\qquad (\xi_1,\xi_2)\in\Z_{L_1}\times\Z_{L_2},$$
for the discrete Fourier transform, and $\theta_A(s)=\sum_x A(x)\overline{A(x+s)}$ for the periodic
autocorrelation. By the Wiener--Khinchin relation $\theta_A$ is the inverse transform of
$|\widehat A|^2$, so $A$ is perfect iff $|\widehat A|$ is constant. For a unimodular array
$\sum_\xi|\widehat A(\xi)|^2=(L_1L_2)^2$, hence
\begin{equation}\label{flat}
A\ \text{unimodular is perfect}\iff|\widehat A(\xi)|=\sqrt{L_1L_2}\ \text{for all }\xi .
\end{equation}

\begin{definition}[type]\label{typedef}
Let $K\ge2$, $\omega=e^{2\pi\sqrt{-1}/K}$, and let $Q$ be an $M_1\times M_2$ array with entries in
$\mu_K$. For $v=(v_1,v_2)\in\Z_K^2$ put
$$\rho_v(x_1,x_2)=e^{-2\pi\sqrt{-1}\,\left(\frac{v_1x_1}{KM_1}+\frac{v_2x_2}{KM_2}\right)},
\qquad 0\le x_i<M_i .$$
$Q$ has \emph{type $v$} if the array $Q\rho_v$ is perfect on $\Z_{M_1}\times\Z_{M_2}$. We write
$\mathcal T_K(M_1,M_2)$ for the set of types $v$ for which an $M_1\times M_2$ array of type $v$ exists.
\end{definition}

Type $(0,0)$ is ordinary perfection. Unwinding the definition, $Q$ has type $v$ iff for every
$s\ne0$
$$\sum_{x}Q(x)\,\overline{Q(x+s)}\;\omega^{-v_1[x_1+s_1\ge M_1]-v_2[x_2+s_2\ge M_2]}=0,$$
i.e.\ the periodic autocorrelation in which every term that wraps around in direction $i$ carries the
weight $\omega^{-v_i}$ vanishes off-peak. For $K=2$ this is the nega-periodic autocorrelation in the
twisted directions, so types $(1,0),(0,1),(1,1)$ are the row-wise quasi-perfect, column-wise
quasi-perfect and doubly quasi-perfect binary arrays of
\cite{JedwabMitchell1988,Wild1988,JMPW1994}, and type $v$ with $v\in\{0,1\}^2$ is Jedwab's
generalised perfect binary array of type $z=v$ \cite[Defs.~2.7, 2.8]{Jedwab1992}: his expansion
$\epsilon(A;z)$ adjoins a negated copy of the array in each direction with $z_i=1$, which is the
twisted extension $E_v$ of Lemma~\ref{twistlem} below, and his arrays are generalised perfect when
the expansion has zero autocorrelation off the subgroup of shifts $u_i\equiv0\pmod{s_i}$. (Jedwab's
``$m$-ary'' arrays for $m>2$ have integer entries in $\{\pm1,\dots\}$ or $\{0,\pm1,\dots\}$, not
roots of unity, so only his binary case is an instance of ours.) For general $K$ and
$v\in\{0,1\}^2$ it is the generalised perfect $K$-phase array of type $z=v$ of
\cite[Defs.~3.1, 3.4]{AEF2024}, whose expansion $\varphi'(g)=\varphi(g\bmod s)+\sum_i\lfloor
g_i/s_i\rfloor$ is again $E_v$. Types with entries outside $\{0,1\}$ do not seem to have been
considered before. A type is an \emph{axis type} if $v_1v_2=0$.

Types only matter when $K$ shares a factor with the side. If $\gcd(K,M_1)=1$, choose $u$ with
$uM_1\equiv1\pmod K$ and write $uM_1=1+Kt$; then
$e^{-2\pi\sqrt{-1}\,v_1x_1/(KM_1)}=\omega^{-uv_1x_1}\,e^{-2\pi\sqrt{-1}\,tv_1x_1/M_1}$, a
$\mu_K$-valued function times a character of $\Z_{M_1}$, and multiplying by a character does not
affect perfection; so $Q$ has type $(v_1,v_2)$ iff $Q\,\omega^{-uv_1x_1}$ has type $(0,v_2)$, and
when both sides are coprime to $K$ every type reduces to perfection. For $K=2$ and odd side this is
\cite[Thm.~7.4]{Jedwab1992} (the modulation $(-1)^{x_1}$), and in general it is the statement that
the carry cocycle is a coboundary when the side is coprime to $K$ \cite[Prop.~4.1(ii)]{AEF2024}.
In the tower of Section~\ref{towersec} every side is a power of $K$, which is the opposite extreme.

There is a group-theoretic reading of types, which we use to compare with the difference-set
literature. A function $f\colon\Z_{M_1}\times\Z_{M_2}\to\Z_K$ with $\omega^f$ of type $v$ is the same thing as a
section of the central extension $0\to\Z_K\to G_v\to\Z_{M_1}\times\Z_{M_2}\to0$ with cocycle given by the
carries $v_1[x_1+y_1\ge M_1]+v_2[x_2+y_2\ge M_2]$ (so $G_{(0,0)}=\Z_{M_1}\times\Z_{M_2}\times\Z_K$ and,
e.g., $G_{(1,0)}=\Z_{KM_1}\times\Z_{M_2}$), whose graph $\{(x,f(x))\}$ has all its ``twisted
differences'' balanced. For $K$ prime a vanishing sum of $M_1M_2$ $K$-th roots of unity is uniform, so
type-$v$ arrays are exactly the semiregular relative $(M_1M_2,K,M_1M_2,M_1M_2/K)$-difference sets in
$G_v$ relative to the central $\Z_K$. For composite $K$ the type condition is weaker than the
difference-set condition. For $K=2$ this is Jedwab's theorem that a generalised perfect binary array
of type $z$ is a relative difference set in an abelian factor group of
$\Z_{(z_1+1)s_1}\times\Z_{(z_2+1)s_2}$ \cite[Thm.~3.2]{Jedwab1992};
the identification of the cocycle with the product of the carry cocycles
$\gamma_m(j,k)=\omega^{\lfloor(j+k)/m\rfloor}$ is in Hughes \cite{Hughes2000} and, for general $K$
and $z\in\{0,1\}^2$, in \cite[Prop.~4.1, Thm.~4.9]{AEF2024}.

\section{Three lemmas}\label{lemmasec}

\begin{lemma}[twisted extension]\label{twistlem}
Let $d_1,d_2,m_1,m_2\ge1$, $N_i=d_im_i$, let $R$ be a complex $m_1\times m_2$ array and
$w\in\Z_{d_1}\times\Z_{d_2}$. Define $E=E^{d}_w[R]$ on $\Z_{N_1}\times\Z_{N_2}$ by
$$E(\bar x_1+m_1y_1,\ \bar x_2+m_2y_2)=e^{2\pi\sqrt{-1}\,(w_1y_1/d_1+w_2y_2/d_2)}\,R(\bar x_1,\bar x_2),
\qquad 0\le\bar x_i<m_i,\ 0\le y_i<d_i .$$
Then $\widehat E(\xi)=0$ unless $\xi_i\equiv w_i\pmod{d_i}$ for $i=1,2$, and
$$\widehat E(w_1+d_1s_1,\,w_2+d_2s_2)=d_1d_2\;\widehat{(R\rho)}(s_1,s_2),\qquad
\rho(\bar x)=e^{-2\pi\sqrt{-1}\,(w_1\bar x_1/N_1+w_2\bar x_2/N_2)},$$
where $\widehat{(R\rho)}$ is the transform on $\Z_{m_1}\times\Z_{m_2}$.
\end{lemma}
\begin{proof}
Substituting $x_i=\bar x_i+m_iy_i$,
$$\widehat E(\xi)=\sum_{\bar x}R(\bar x)\,e^{-2\pi\sqrt{-1}\sum_i\xi_i\bar x_i/N_i}
\prod_{i=1,2}\ \sum_{y_i=0}^{d_i-1}e^{2\pi\sqrt{-1}\,(w_i-\xi_i)y_i/d_i}.$$
The inner geometric sums equal $d_i[\xi_i\equiv w_i\ (d_i)]$. For $\xi_i=w_i+d_is_i$ the exponent
$\xi_i\bar x_i/N_i=w_i\bar x_i/N_i+s_i\bar x_i/m_i$ splits into the factor $\rho$ and the
$\Z_{m_1}\times\Z_{m_2}$ character $s$.
\end{proof}

With $d_1=d_2=K$ and $w=v$ the factor $\rho$ is exactly $\rho_v$ of Definition~\ref{typedef}. Hence,
writing $E_v[Q]=E^{(K,K)}_v[Q]$ and calling an array on $\Z_{N_1}\times\Z_{N_2}$ \emph{coset-flat}
for the coset $v+(K\Z_{N_1}\times K\Z_{N_2})$ if its transform vanishes off that coset and has constant
modulus on it:

\begin{corollary}\label{cosetflat}
$Q$ has type $v$ if and only if $E_v[Q]$ is coset-flat for the coset $v+K\Z^2$.
\end{corollary}

\begin{lemma}[interleaving]\label{interlem}
Let $F^{(a)}$, $0\le a<K$, be complex arrays on $\Z_{N_1}\times\Z_{N_2}$ and define $S$ on
$\Z_{KN_1}\times\Z_{N_2}$ by $S(Ki'+a,x_2)=F^{(a)}(i',x_2)$. Then
$$\widehat S(\xi_1,\xi_2)=\sum_{a=0}^{K-1}e^{-2\pi\sqrt{-1}\,\xi_1a/(KN_1)}\,
\widehat{F^{(a)}}(\xi_1\bmod N_1,\ \xi_2).$$
Consequently, if the supports of the $\widehat{F^{(a)}}$ are pairwise disjoint and each
$\widehat{F^{(a)}}$ has constant modulus $C$ on its support, then $|\widehat S(\xi)|\in\{0,C\}$, and
$|\widehat S(\xi)|=C$ exactly when $(\xi_1\bmod N_1,\xi_2)$ lies in the union of the supports.
\end{lemma}
\begin{proof}
$\sum_{i',a}F^{(a)}(i',x_2)e^{-2\pi\sqrt{-1}(\xi_1(Ki'+a)/(KN_1)+\xi_2x_2/N_2)}$; the term $Ki'$
contributes the $\Z_{N_1}$-character $\xi_1\bmod N_1$ and the term $a$ the twiddle factor.
\end{proof}

Interleaving in both directions is the composite of two one-directional interleavings. In
particular, if $K^2$ arrays $F^{(a,b)}$ on $\Z_N^2$ are coset-flat for the $K^2$ distinct cosets of
$K\Z_N^2$ with a common modulus, then $S(Ki'+a,Kj'+b)=F^{(a,b)}(i',j')$ has constant-modulus
transform on $\Z_{KN}^2$, hence is perfect by \eqref{flat}. Combined with
Corollary~\ref{cosetflat} this already gives a perfect array of side $K^2M$ from arrays of every type
at side $M$ (Construction~\ref{thmA} below); the typed lift of the next section is the sharper,
one-directional statement needed to iterate.

Lemma~\ref{interlem} contains the Generalised Array Orthogonality Property (GAOP) of
\cite{Blake2017} as a special case. Two arrays are orthogonal under all cyclic shifts exactly when
the pointwise product of their transforms vanishes, i.e.\ when their spectral supports are disjoint,
and a family of arrays is periodic-complementary exactly when the sum of the squared moduli of their
transforms is constant; so sub-arrays that are pairwise orthogonal and complementary have disjoint
flat supports covering the frequency plane, and Lemma~\ref{interlem} says the interleaved array is
perfect. This is Theorem~7.0.3 of \cite{Blake2017}. The lemma is stronger in that the same
decomposition also certifies arrays of nonzero type, which the GAOP, whose conclusion is always
perfection, cannot do; this is what makes the induction of Section~\ref{towersec} possible.

\begin{lemma}[shear with carry]\label{shearlem}
Let $Q$ be an $M\times M$ array over $\mu_K$ of type $v=(v_1,v_2)$. Then
\begin{equation}\label{shear1}
Q'(x,y)=\omega^{\,v_2\fl{(x+y)/M}}\;Q\big(x,\ (x+y)\bmod M\big)
\end{equation}
has type $(v_1+v_2,\ v_2)$, and the transposed shear
$Q''(x,y)=\omega^{\,v_1\lfloor(x+y)/M\rfloor}Q((x+y)\bmod M,\ y)$ has type $(v_1,\ v_1+v_2)$.
\end{lemma}
\begin{proof}
Let $E=E_v[Q]$ on $\Z_{KM}^2$ and $E'(x,y)=E(x,x+y)$, the image of $E$ under the automorphism
$(x,y)\mapsto(x,x+y)$ of $\Z_{KM}^2$. Then $\widehat{E'}(\xi_1,\xi_2)=\widehat E(\xi_1-\xi_2,\xi_2)$,
so $E'$ is coset-flat for the coset $\{\xi:\xi_1-\xi_2\equiv v_1,\ \xi_2\equiv v_2\}=
(v_1+v_2,v_2)+K\Z^2$. On the other hand, by the twisted periodicity of $E$ in both directions,
$E'(x+M,y)=\omega^{v_1+v_2}E'(x,y)$ and $E'(x,y+M)=\omega^{v_2}E'(x,y)$, so $E'=E_{(v_1+v_2,v_2)}[Q']$
where $Q'$ is the restriction of $E'$ to $[0,M)^2$, namely
$Q'(x,y)=E(x,x+y)=\omega^{v_2\lfloor(x+y)/M\rfloor}Q(x,(x+y)\bmod M)$. Corollary~\ref{cosetflat}
gives the type of $Q'$. The transposed statement follows by symmetry.
\end{proof}

For $K=2$ the shear is Jedwab's transform $\tau(A;c)$ \cite[Def.~5.7]{Jedwab1992},
$d[l,j]=(-1)^{y}a[l-cj,\,j]$ with $y=1$ exactly when the index wraps, and
\cite[Thm.~7.6]{Jedwab1992} states that for $ct\equiv s\pmod{2s}$ it carries type $(1,0)$ to type
$(1,1)$ on an $s\times t$ array (so $c=1$ for square arrays); this generalises Wild's theorem that
quasi-perfect and doubly quasi-perfect binary arrays coexist when $t/\gcd(s,t)$ is odd
\cite{Wild1988}, \cite[Result~4.5]{Jedwab1992}.

The two shears are the elementary matrices $\begin{pmatrix}1&1\\0&1\end{pmatrix}$ and
$\begin{pmatrix}1&0\\1&1\end{pmatrix}$ acting on the row vector $v$; they generate $SL_2(\Z)$, and the
image of $SL_2(\Z)$ in $SL_2(\Z_K)$ is all of $SL_2(\Z_K)$.

\begin{lemma}[orbits]\label{orbitlem}
The orbits of $SL_2(\Z_K)$ on $\Z_K^2$ are the sets $\{(x,y):\gcd(x,y,K)=g\}$, $g\mid K$. Hence every
type is reachable from an axis type $(g,0)$ by a finite sequence of the shears of
Lemma~\ref{shearlem}.
\end{lemma}
\begin{proof}
$\gcd(x,y,K)$ is invariant. For the converse reduce by the Chinese remainder theorem to $K=p^e$:
a vector with $\min(\nu_p(x),\nu_p(y))=j$ is $p^j$ times a primitive vector, and $SL_2(\Z_{p^e})$ is
transitive on primitive vectors since any primitive vector is the first row of an invertible matrix,
whose determinant can be normalised in the second row.
\end{proof}

Before using these lemmas we record why the obvious recursion fails. If $A$ is a perfect array on a group $G$ of order $|G|$, its $|G|$ cyclic shifts are an orthogonal
basis of $\C^{G}$. Consequently no nonzero array is orthogonal to all shifts of $A$: two perfect arrays can never
have disjoint spectral supports. Components with disjoint supports, as Lemma~\ref{interlem}
requires, therefore cannot be perfect; they can only be complementary. This rules out any recursion
that passes from perfect components to a perfect array, and explains why the Frank--Heimiller array $\omega^{\lfloor x_1x_2/K\rfloor}$ of side $K^2$, whose
$K$-polyphase components are constants times a Frank--Heimiller core, cannot be enlarged by the same
device. In the construction below the components are the twisted extensions $E_v[Q]$ of
Lemma~\ref{twistlem}: their spectra live on disjoint cosets, but none of them is perfect.

Conversely, one checks that an array obtained by interleaving $K^2$ coset-flat arrays as above never
has a nonzero type: a nonzero type of the interleaved array would force the tiles to be twisted by
$K^2$-th roots of unity. This is why the recursion of the next section is formulated for a
one-directional lift of the \emph{twisted extension} $E_v[S]$ rather than of $S$ itself.

\section{The lift and the tower}\label{towersec}

Fix $K\ge2$. All index functions below are integer valued and read modulo $K$; the array is
$\omega^{f}$.

\begin{theorem}[one-directional lift]\label{liftthm}
Let $v_1\in\Z_K$ and, for each $b\in\Z_K$, let $R^{(b)}$ be an $m_1\times m_2$ array over $\mu_K$ of
type $(v_1,b)$. Define the $(Km_1)\times(Km_2)$ array
\begin{equation}\label{liftdef}
S\big(Ki'+a,\ \bar x_2+m_2y_2\big)=\omega^{\,a\,y_2}\,R^{(a)}(i',\bar x_2),
\end{equation}
where $0\le i'<m_1$, $0\le a<K$, $0\le\bar x_2<m_2$ and $0\le y_2<K$.
Then $S$ has type $(v_1,0)$. Transposing, from arrays $R^{(a)}$ of type $(a,v_2)$ one obtains an
array of type $(0,v_2)$.
\end{theorem}
\begin{proof}
Let $E=E_{(v_1,0)}[S]$ on $\Z_{K^2m_1}\times\Z_{K^2m_2}$, so that
$E(x_1+Km_1y_1',\ x_2+Km_2y_2')=\omega^{v_1y_1'}S(x_1,x_2)$. Decimate $E$ in the first direction:
$E^{(a)}(i,x_2)=E(Ki+a,x_2)$ for $i\in\Z_{Km_1}$, $x_2\in\Z_{K^2m_2}$. Write $i=i'+m_1y_1$ with
$0\le i'<m_1$, $0\le y_1<K$, and $x_2=\bar x_2+m_2y_2$ with $0\le\bar x_2<m_2$, $0\le y_2<K^2$. Since
$Ki+a=(Ki'+a)+Km_1y_1$, the twisted periodicity of $E$ and \eqref{liftdef} give
$$E^{(a)}(i'+m_1y_1,\ \bar x_2+m_2y_2)=\omega^{\,v_1y_1}\,\omega^{\,a(y_2\bmod K)}R^{(a)}(i',\bar x_2)
=e^{2\pi\sqrt{-1}\,\left(\frac{v_1y_1}{K}+\frac{Ka\,y_2}{K^2}\right)}R^{(a)}(i',\bar x_2),$$
i.e.\ $E^{(a)}=E^{(K,K^2)}_{(v_1,Ka)}[R^{(a)}]$ in the notation of Lemma~\ref{twistlem}. By that
lemma $\widehat{E^{(a)}}$ is supported on $\{\xi_1\equiv v_1\ (K),\ \xi_2\equiv Ka\ (K^2)\}\subset
\Z_{Km_1}\times\Z_{K^2m_2}$, where it equals $K^3$ times the transform of $R^{(a)}\rho$ with
$\rho(\bar x)=e^{-2\pi\sqrt{-1}(v_1\bar x_1/(Km_1)+Ka\bar x_2/(K^2m_2))}=\rho_{(v_1,a)}(\bar x)$; this
has constant modulus because $R^{(a)}$ has type $(v_1,a)$. The $K$ supports are pairwise disjoint and
their union is $\{\xi_1\equiv v_1\ (K),\ \xi_2\equiv0\ (K)\}$. By Lemma~\ref{interlem},
$\widehat E$ vanishes off the coset $(v_1,0)+K\Z^2$ of $\Z_{K^2m_1}\times\Z_{K^2m_2}$ and has
constant modulus on it: $E$ is coset-flat, so $S$ has type $(v_1,0)$ by Corollary~\ref{cosetflat}.
\end{proof}

In terms of index functions, with $f^{(v)}_{m}$ denoting an exponent of type $v$ at shape $m$,
\eqref{liftdef} reads
\begin{equation}\label{lift1}
f^{(v_1,0)}_{Km}(x_1,x_2)=(x_1\bmod K)\fl{\frac{x_2}{m_2}}+f^{(v_1,\ x_1\bmod K)}_{m}
\Big(\fl{\frac{x_1}{K}},\ x_2\bmod m_2\Big),
\end{equation}
and the transposed lift is
\begin{equation}\label{lift2}
f^{(0,v_2)}_{Km}(x_1,x_2)=(x_2\bmod K)\fl{\frac{x_1}{m_1}}+f^{(x_2\bmod K,\ v_2)}_{m}
\Big(x_1\bmod m_1,\ \fl{\frac{x_2}{K}}\Big).
\end{equation}

For $K=2$ the lift \eqref{liftdef} with $(v_1,0)=(0,0)$ reads: even columns of $S$ are the
perfect array $R^{(0)}$, odd columns are the row-wise quasi-perfect array $R^{(1)}$ in the upper
half and $-R^{(1)}$ in the lower half. This is the doubling PBA$(s,t)$ + QPBA$(s,t)$ $\Rightarrow$
PBA$(2s,2t)$ of Jedwab and Mitchell \cite{JedwabMitchell1988}, \cite[Result~4.3]{Jedwab1992},
and the lift with $v_1=1$ is their QPBA$(s,t)$ + DQPBA$(s,t)$ $\Rightarrow$ QPBA$(2s,2t)$
\cite[Thm.~1]{JedwabMitchell1990}, \cite[Result~4.4]{Jedwab1992}, whose proof is exactly ``interleave
the rows of $[B\ B]$ and $[C\ {-C}]$''. Both are special cases of \cite[Thm.~7.5]{Jedwab1992}, which
is Theorem~\ref{liftthm} for $K=2$ in any dimension (with the twist in the first coordinate) and is
moreover an equivalence: the interleaved array has the output type if and only if the two inputs have
the input types. Theorem~\ref{liftthm} is thus the $K$-ary version of the Jedwab--Mitchell doubling,
with the two binary inputs replaced by $K$ inputs of types $(v_1,0),\dots,(v_1,K-1)$ and the sign
$-1$ by the twists $\omega^{ay_2}$; we have not pursued the converse.

\begin{lemma}[base]\label{baselem}
The Frank array $F(x,y)=\omega^{xy}$ on $\Z_K^2$ \cite{Heimiller1961,Frank1962} has every axis type $(a,0)$ and $(0,b)$,
$a,b\in\Z_K$. Hence, by Lemma~\ref{orbitlem}, $\mathcal T_K(K,K)=\Z_K^2$.
\end{lemma}
\begin{proof}
With $\zeta=e^{2\pi\sqrt{-1}/K^2}$, $F\rho_{(a,0)}(x,y)=\omega^{xy}\zeta^{-ax}$. Summing first over
$y$, $\sum_y\omega^{xy}\omega^{-\beta y}=K[x=\beta]$, so
$\widehat{(F\rho_{(a,0)})}(\alpha,\beta)=K\zeta^{-a\beta}\omega^{-\alpha\beta}$, of modulus $K$ for
all $(\alpha,\beta)$. The type $(0,b)$ is the transposed statement.
\end{proof}

\begin{theorem}[tower]\label{towerthm}
If $\mathcal T_K(M,M)=\Z_K^2$ then $\mathcal T_K(KM,KM)=\Z_K^2$. Consequently
$\mathcal T_K(K^n,K^n)=\Z_K^2$ for every $n\ge1$.
\end{theorem}
\begin{proof}
Given arrays of every type at side $M$, Theorem~\ref{liftthm} produces, for each $v_1$, an array of
type $(v_1,0)$ at side $KM$ from the arrays of types $(v_1,b)$, $b\in\Z_K$; transposed, it produces
every type $(0,v_2)$. In particular every axis type $(g,0)$ exists at side $KM$, and
Lemma~\ref{orbitlem} with Lemma~\ref{shearlem} supplies the remaining types. Induction from
Lemma~\ref{baselem} gives all $n$.
\end{proof}

The type-$0$ member of the family is the result announced in the introduction.

\begin{theorem}\label{mainthm}
For every integer $K\ge2$ and every $n\ge1$ there is a perfect $K^n\times K^n$ array over the $K$-th
roots of unity, given explicitly by the recursion \eqref{lift1}--\eqref{shear1}.
\end{theorem}
\begin{proof}
Take $v=(0,0)$ in Theorem~\ref{towerthm}: type $(0,0)$ is perfection by Definition~\ref{typedef}.
The recursion \eqref{lift1}--\eqref{shear1} is the proof of Theorem~\ref{towerthm} made explicit;
this is Construction~\ref{tower} below, where the claim is proved again in that form.
\end{proof}

We emphasise again that the existence part of Theorem~\ref{towerthm} is known. For $K=p$ prime,
the group-theoretic reading of Section~\ref{defsec} identifies the arrays of Theorem~\ref{towerthm} with
$(p^{2n},p,p^{2n},p^{2n-1})$-relative difference sets in $\Z_{p^n}^2\times\Z_p$ (type $0$) and in
$\Z_{p^{n+1}}\times\Z_{p^n}$ (the other types); both groups have exponent at most $p^{n+1}$, so
their existence, for every choice of the forbidden subgroup of order $p$, is the case $c=n$ of
\cite[Thm.~2.4]{MaSchmidt1995}; the type-$0$ groups fall under their Theorem~2.2 (the $K$-matrix
construction of Davis and Kraemer) and the groups $\Z_{p^{n+1}}\times\Z_{p^n}$ with forbidden
subgroup $\langle p^n\rangle$ under their Theorem~2.1 (Davis's construction, with
$E=\langle(p,0)\rangle\times\langle(p,1)\rangle$). The split between the two cases is exact: Davis's
original theorem \cite[Thm.~4.1]{Davis1992} (order $p^{j+2}$, $j$ odd, exponent at most
$p^{(j+3)/2}$) already gives every nonzero type directly, with $j=2n-1$, but its proof places the
forbidden subgroup inside a cyclic factor of largest order, as Davis remarks, so it does not reach
the type-$0$ group, where the forbidden subgroup is the direct factor $\Z_p$; that case needs the
$K$-matrix method. For $K=2$ every type at every side $2^n$ is \cite[Thm.~9.2(ii)]{Jedwab1992},
and type $0$ alone is \cite[Thm.~4.2]{Davis1991}. An explicit $p$-ary construction on $\Z_{p^n}^2$ is also given by
Do Duc \cite[Thm.~4.3]{DoDuc2021} through the local ring $\Z_{p^n}$ (the sets
$\{(x,y):f_i(\varphi(x)y)=1\}$ with $\varphi(p^ku)=p^ku^{-1}$, in the style of Dillon's binary formula
\cite{Dillon1990}); the nonzero axis types in $\Z_{p^{n+1}}\times\Z_{p^n}$ are also produced by
the local-ring construction of Ma and Schmidt \cite[Thm.~4.2, Cor.~4.3]{MaSchmidt2000} applied to
the chain ring $\Z_{p^{n+1}}$, in the generality of Hou and Sehgal \cite[Thm.~3.2]{HouSehgal2000}. For $K=p^e$ a $p$-ary perfect array is an array over $K$-th roots,
and for $K=\prod p_i^{e_i}$ the arrays for the prime-power parts combine, since
$\Z_{K^n}\cong\prod\Z_{p_i^{e_in}}$, into a perfect $K^n\times K^n$ array over $K$-th roots by
L\"uke's product construction \cite{Luke1988}. The content of Theorem~\ref{towerthm} is therefore the
uniform formula and the fact that the alphabet is exactly $\mu_K$ at every level, not the existence.
Note that the bilinear-form construction of \cite{DoDucSchmidt2019} gives $K$-ary arrays on
$\Z_{K^n}^2$ only for $n\le2$ (it requires $\nu_p(K)\ge\lceil n\nu_p(K)/2\rceil$).

\begin{construction}\label{tower}
\emph{(The explicit recursion.)} Let $f^{(a,0)}_K=f^{(0,b)}_K=xy$ for all $a,b\in\Z_K$ and complete
the set $\{f^{(v)}_K\}$ to all $v\in\Z_K^2$ by the shears \eqref{shear1} (and its transpose),
following any path provided by Lemma~\ref{orbitlem}. Given $\{f^{(v)}_M\}_{v\in\Z_K^2}$, form
$f^{(v_1,0)}_{KM}$ by \eqref{lift1}, $f^{(0,v_2)}_{KM}$ by \eqref{lift2}, and complete by shears. Then
$\omega^{f^{(0,0)}_{K^n}}$ is a perfect $K^n\times K^n$ array over $\mu_K$.
\end{construction}
\begin{proof}
We show by induction on $n$ that $\omega^{f^{(v)}_{K^n}}$ has type $v$ for every $v\in\Z_K^2$; the
claim is the case $v=(0,0)$, which is perfection by Definition~\ref{typedef}. For $n=1$ the axis
types are Lemma~\ref{baselem}, and Lemma~\ref{shearlem} shows that each shear carries an array of
type $(v_1,v_2)$ to one of type $(v_1+v_2,v_2)$ (resp.\ $(v_1,v_1+v_2)$), so the array reached from
an axis type by a sequence of shears has the type that Lemma~\ref{orbitlem} assigns to it, and the
completion reaches every $v$. Suppose the claim holds at side $M=K^{n-1}$. Writing $x_1=Ki'+a$ and
$x_2=\bar x_2+m_2y_2$, so that $a=x_1\bmod K$, $i'=\lfloor x_1/K\rfloor$, $y_2=\lfloor x_2/m_2\rfloor$
and $\bar x_2=x_2\bmod m_2$, \eqref{lift1} is \eqref{liftdef} written in exponents with
$R^{(a)}=\omega^{f^{(v_1,a)}_M}$; by the induction hypothesis $R^{(a)}$ has type $(v_1,a)$, so
Theorem~\ref{liftthm} gives $\omega^{f^{(v_1,0)}_{KM}}$ type $(v_1,0)$. Likewise \eqref{lift2} is the
transposed lift and gives type $(0,v_2)$. The completing shears then produce every remaining type at
side $KM$ exactly as at $n=1$. Note that Theorem~\ref{liftthm} and Lemma~\ref{shearlem} ask only
that their inputs have the stated types, so the result does not depend on the shear paths chosen in
the completion steps: different paths give different arrays, all of the asserted types.
\end{proof}

The completion step is a breadth-first search over the two shears, which terminates by
Lemma~\ref{orbitlem}.

\begin{construction}\label{thmA}
\emph{(Two-step version.)} From arrays $Q^{(u)}$ of every type $u\in\Z_K^2$ at side $M$ the array
$S(Ki'+a,Kj'+b)=E_{(a,b)}[Q^{(a,b)}](i',j')$ is a perfect array of side $K^2M$: its $K^2$ polyphase
components are coset-flat for the $K^2$ cosets of $K\Z_{KM}^2$ (Corollary~\ref{cosetflat}), each
with the same modulus $K^2M$ by Lemma~\ref{twistlem} and \eqref{flat}, and the two-directional
form of Lemma~\ref{interlem} applies. Two successive lifts of Construction~\ref{tower} (one of
each direction, in either order) give the same array up to a permutation of the polyphase
components: in the lift \eqref{liftdef} the residue $a$ of the interleaved coordinate carries the
twist in the \emph{other} coordinate and the \emph{other} entry of the type, so the composite places
$E_{(b,a)}[Q^{(b,a)}]$, not $E_{(a,b)}[Q^{(a,b)}]$, at the residue $(a,b)$. Both arrangements are
perfect by the same argument, since $(a,b)\mapsto(b,a)$ permutes the cosets.
\end{construction}

\begin{figure}[H]
\centering
\begin{tikzpicture}[x=1cm,y=1cm,>=stealth,
  cell/.style={draw,minimum size=0.95cm,inner sep=0pt,font=\footnotesize},
  lifted/.style={cell,fill=black!25},
  sheared/.style={cell,fill=black!8},
  lab/.style={font=\small}]
\begin{scope}[shift={(0,0)}]
\node[lab] at (1.95,4.2) {(a) one lift, $K=3$};
\foreach \a/\y in {0/2,1/1,2/0}{
  \foreach \b/\x in {0/0,1/1,2/2}{
    \pgfmathtruncatemacro{\e}{mod(\a*\b,3)}
    \node[cell,minimum size=1.3cm] at (\x*1.3+0.65,\y*1.3+0.65) {$\omega^{\e}R^{(\a)}$};
  }
  \node[lab,anchor=east] at (-0.1,\y*1.3+0.65) {$x_1\equiv \a$};
}
\foreach \b/\x in {0/0,1/1,2/2}{\node[lab] at (\x*1.3+0.65,-0.3) {$y_2=\b$};}
\node[lab,anchor=north] at (1.95,-0.6) {rows gathered by $a=x_1\bmod K$};
\end{scope}
\begin{scope}[shift={(5.6,0.3)}]
\node[lab] at (1.6,4.15) {(b) side $M$: all types};
\foreach \u in {0,1,2}{\foreach \v in {0,1,2}{
  \node[cell] (L\u\v) at (\v*1.1+0.5,{(2-\u)*1.1+0.5}) {$(\u,\v)$};}}
\end{scope}
\begin{scope}[shift={(10.4,0.3)}]
\node[lab] at (1.6,4.15) {side $KM$: all types};
\foreach \u in {0,1,2}{\foreach \v in {0,1,2}{
  \ifnum\u=0 \node[lifted] (R\u\v) at (\v*1.1+0.5,{(2-\u)*1.1+0.5}) {$(\u,\v)$};
  \else\ifnum\v=0 \node[lifted] (R\u\v) at (\v*1.1+0.5,{(2-\u)*1.1+0.5}) {$(\u,\v)$};
  \else \node[sheared] (R\u\v) at (\v*1.1+0.5,{(2-\u)*1.1+0.5}) {$(\u,\v)$};\fi\fi }}
\foreach \u in {1,2}{\draw[->,thick] (R\u0) -- (R\u1); \draw[->,thick] (R\u1) -- (R\u2);}
\end{scope}
\foreach \u in {0,1,2}{
  \draw[->,thick,black!60] (9.0,{(2-\u)*1.1+0.5+0.3}) -- (R\u0.west);}
\foreach \v in {1,2}{
  \draw[->,thick,black!60] ({5.6+\v*1.1+0.5},3.55) to[out=90,in=90,looseness=0.5] (R0\v.north);}
\node[lab,anchor=north,align=center] at (9.8,-0.05)
  {rows $\{(v_1,b)\}_{b}\to(v_1,0)$ by \eqref{lift1};\ columns $\to(0,v_2)$ by \eqref{lift2}\\
   dark: from lifts;\ light: from the shears \eqref{shear1}};
\end{tikzpicture}
\caption{The recursion. (a) The lift \eqref{liftdef} for $K=3$: the rows of $S$ with
$x_1\equiv a\pmod K$ (interleaved in $S$, gathered here) are the array $R^{(a)}$ of type $(v_1,a)$
tiled $K$ times along $x_2$ with the twist $\omega^{ay_2}$ on the $y_2$-th tile; the result has type
$(v_1,0)$. (b) The type bookkeeping for one level ($v_1$ down, $v_2$ across): each row of types
$(v_1,b)$, $b\in\Z_K$, at side $M$ lifts to the axis type $(v_1,0)$ at side $KM$, each column to
$(0,v_2)$, and the shears fill the remaining types from the axis types.}\label{figrec}
\end{figure}
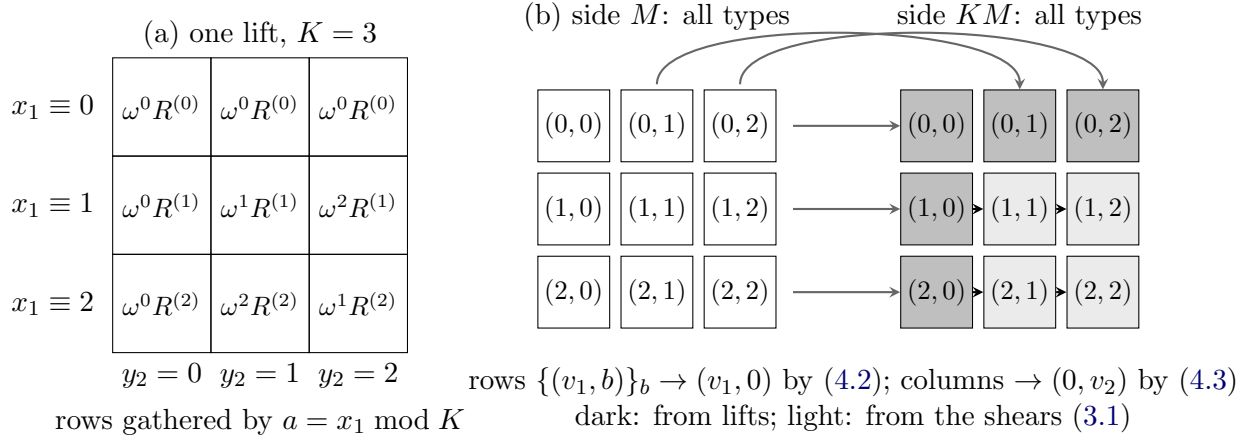

Figure~\ref{figrec} shows one lift and the type bookkeeping of one level of the recursion.

\subsection*{Why the tower is two-dimensional}\label{onedsec}

Nothing in Lemmas~\ref{twistlem} and \ref{interlem} needs two coordinates, so it is worth seeing
exactly where the same construction fails for sequences, since perfect sequences of unbounded
length at a fixed alphabet do not exist (Section~\ref{litsec}). A type is now a single residue
$b\in\Z_K$: a sequence $Q$ of length $M$ over $\mu_K$ has type $b$ if $Q(x)e^{-2\pi\sqrt{-1}\,bx/(KM)}$
is perfect, i.e.\ if its autocorrelation with the wrapped terms weighted by $\omega^{-b}$ vanishes
off-peak. The one-dimensional lift is forced to interleave and twist the same coordinate: from
sequences $R^{(b)}$ of length $M$ and type $b$, $b\in\Z_K$, form the twisted extensions
$T^{(b)}(\bar x+My)=\omega^{by}R^{(b)}(\bar x)$ of length $KM$ and interleave them,
$S(Ki'+a)=T^{(a)}(i')$. Exactly as in Theorem~\ref{liftthm}, the $T^{(b)}$ have flat spectra on the
$K$ cosets of $K\Z_{KM}$, so $S$ is a perfect sequence of length $K^2M$. Starting from $M=1$, where
the one-term sequence $[1]$ has every type, this gives $S(Ky+a)=\omega^{ay}$: the Frank sequence of
length $K^2$. (The two-dimensional analogue from side $1$ is Construction~\ref{thmA} with $M=1$,
which gives $\omega^{ay_1+by_2}$, the Frank--Heimiller array of side $K^2$.)

The recursion stops there, for a structural reason. The lift zeroes the type in the twisted
coordinate, and in one dimension that is the only coordinate: the one-dimensional lift produces
perfect sequences and nothing else, so its output can never be fed back in. In two dimensions the
lift keeps the type of the interleaved coordinate (Theorem~\ref{liftthm} outputs $(v_1,0)$ with
$v_1$ arbitrary), and the shears of Lemma~\ref{shearlem}, which need a second coordinate to mix
with, regenerate every nonzero type from the axis types; there is no $SL_1$. So the sequence tower
would need the nonzero types to be supplied from outside at each level, and they are not available:
a type-$b$ sequence of length $M$ over $\mu_K$ is, for $K=p$ prime and $b\ne0$, a
$(M,p,M,M/p)$-relative difference set in the cyclic group $\Z_{pM}$ relative to its subgroup of
order $p$ (the group-theoretic reading of Section~\ref{defsec} with the carry cocycle in one
coordinate). At $M=p^2$, the length of the Frank sequence, such a relative difference set would live
in $\Z_{p^3}$, whose exponent $p^3$ exceeds the bound $p^{c+1}=p^2$ of
\cite[Thm.~2.4]{MaSchmidt1995} ($c=1$); so no perfect sequence of length $p^2$ over $\mu_p$ has a
nonzero type, and not even one further lift is possible. (At length $p$ the exponent bound does not
decide, but for $p=3$ one checks directly that no sequence of length $3$ over the cube roots of
unity has a nonzero type.) The two-dimensional tower exists because the twisted directions can be
rotated into one another; a sequence has nowhere to put the twist.

\section{Examples}

\begin{example}[side $K^2$]
Applying \eqref{lift2} to the base gives
$$f^{(0,0)}_{K^2}(x_1,\ Kj'+b)=b\fl{\frac{x_1}{K}}+(x_1\bmod K)\,j' .$$
Comparing with the Frank--Heimiller form $\lfloor x_1x_2/K\rfloor=x_1j'+\lfloor x_1b/K\rfloor$ one finds
$f^{(0,0)}_{K^2}(x_1,x_2)\equiv\lfloor x_1x_2/K\rfloor-\lfloor(x_1\bmod K)(x_2\bmod K)/K\rfloor\pmod K$:
the Frank--Heimiller array with the residue term $\lfloor(x_1\bmod K)(x_2\bmod K)/K\rfloor$ removed.
That term is inert: it depends only on the residues $(x_1\bmod K,x_2\bmod K)$, i.e.\ it is constant
on each $K$-polyphase component, and multiplying a component by a constant changes neither the
support nor the modulus of its spectrum.
\end{example}

\begin{figure}[H]
\centering
\begin{tikzpicture}[x=0.28cm,y=0.28cm]
\fill (4,-1) rectangle (5,-2);
\fill (5,-1) rectangle (6,-2);
\fill (6,-1) rectangle (7,-2);
\fill (7,-1) rectangle (8,-2);
\fill (8,-1) rectangle (9,-2);
\fill (9,-1) rectangle (10,-2);
\fill (10,-1) rectangle (11,-2);
\fill (11,-1) rectangle (12,-2);
\fill (4,-2) rectangle (5,-3);
\fill (5,-2) rectangle (6,-3);
\fill (6,-2) rectangle (7,-3);
\fill (7,-2) rectangle (8,-3);
\fill (12,-2) rectangle (13,-3);
\fill (13,-2) rectangle (14,-3);
\fill (14,-2) rectangle (15,-3);
\fill (15,-2) rectangle (16,-3);
\fill (3,-3) rectangle (4,-4);
\fill (6,-3) rectangle (7,-4);
\fill (7,-3) rectangle (8,-4);
\fill (8,-3) rectangle (9,-4);
\fill (9,-3) rectangle (10,-4);
\fill (10,-3) rectangle (11,-4);
\fill (12,-3) rectangle (13,-4);
\fill (13,-3) rectangle (14,-4);
\fill (2,-4) rectangle (3,-5);
\fill (3,-4) rectangle (4,-5);
\fill (6,-4) rectangle (7,-5);
\fill (7,-4) rectangle (8,-5);
\fill (10,-4) rectangle (11,-5);
\fill (11,-4) rectangle (12,-5);
\fill (14,-4) rectangle (15,-5);
\fill (15,-4) rectangle (16,-5);
\fill (2,-5) rectangle (3,-6);
\fill (3,-5) rectangle (4,-6);
\fill (4,-5) rectangle (5,-6);
\fill (7,-5) rectangle (8,-6);
\fill (8,-5) rectangle (9,-6);
\fill (9,-5) rectangle (10,-6);
\fill (13,-5) rectangle (14,-6);
\fill (14,-5) rectangle (15,-6);
\fill (2,-6) rectangle (3,-7);
\fill (4,-6) rectangle (5,-7);
\fill (5,-6) rectangle (6,-7);
\fill (7,-6) rectangle (8,-7);
\fill (10,-6) rectangle (11,-7);
\fill (12,-6) rectangle (13,-7);
\fill (13,-6) rectangle (14,-7);
\fill (15,-6) rectangle (16,-7);
\fill (2,-7) rectangle (3,-8);
\fill (5,-7) rectangle (6,-8);
\fill (7,-7) rectangle (8,-8);
\fill (8,-7) rectangle (9,-8);
\fill (9,-7) rectangle (10,-8);
\fill (11,-7) rectangle (12,-8);
\fill (12,-7) rectangle (13,-8);
\fill (14,-7) rectangle (15,-8);
\fill (1,-8) rectangle (2,-9);
\fill (3,-8) rectangle (4,-9);
\fill (5,-8) rectangle (6,-9);
\fill (7,-8) rectangle (8,-9);
\fill (9,-8) rectangle (10,-9);
\fill (11,-8) rectangle (12,-9);
\fill (13,-8) rectangle (14,-9);
\fill (15,-8) rectangle (16,-9);
\fill (1,-9) rectangle (2,-10);
\fill (3,-9) rectangle (4,-10);
\fill (4,-9) rectangle (5,-10);
\fill (6,-9) rectangle (7,-10);
\fill (8,-9) rectangle (9,-10);
\fill (10,-9) rectangle (11,-10);
\fill (13,-9) rectangle (14,-10);
\fill (15,-9) rectangle (16,-10);
\fill (1,-10) rectangle (2,-11);
\fill (3,-10) rectangle (4,-11);
\fill (4,-10) rectangle (5,-11);
\fill (6,-10) rectangle (7,-11);
\fill (9,-10) rectangle (10,-11);
\fill (11,-10) rectangle (12,-11);
\fill (12,-10) rectangle (13,-11);
\fill (14,-10) rectangle (15,-11);
\fill (1,-11) rectangle (2,-12);
\fill (5,-11) rectangle (6,-12);
\fill (6,-11) rectangle (7,-12);
\fill (8,-11) rectangle (9,-12);
\fill (10,-11) rectangle (11,-12);
\fill (11,-11) rectangle (12,-12);
\fill (12,-11) rectangle (13,-12);
\fill (15,-11) rectangle (16,-12);
\fill (1,-12) rectangle (2,-13);
\fill (2,-12) rectangle (3,-13);
\fill (5,-12) rectangle (6,-13);
\fill (6,-12) rectangle (7,-13);
\fill (9,-12) rectangle (10,-13);
\fill (10,-12) rectangle (11,-13);
\fill (13,-12) rectangle (14,-13);
\fill (14,-12) rectangle (15,-13);
\fill (1,-13) rectangle (2,-14);
\fill (2,-13) rectangle (3,-14);
\fill (4,-13) rectangle (5,-14);
\fill (5,-13) rectangle (6,-14);
\fill (8,-13) rectangle (9,-14);
\fill (11,-13) rectangle (12,-14);
\fill (14,-13) rectangle (15,-14);
\fill (15,-13) rectangle (16,-14);
\fill (1,-14) rectangle (2,-15);
\fill (2,-14) rectangle (3,-15);
\fill (3,-14) rectangle (4,-15);
\fill (4,-14) rectangle (5,-15);
\fill (9,-14) rectangle (10,-15);
\fill (10,-14) rectangle (11,-15);
\fill (11,-14) rectangle (12,-15);
\fill (12,-14) rectangle (13,-15);
\fill (1,-15) rectangle (2,-16);
\fill (2,-15) rectangle (3,-16);
\fill (3,-15) rectangle (4,-16);
\fill (8,-15) rectangle (9,-16);
\fill (12,-15) rectangle (13,-16);
\fill (13,-15) rectangle (14,-16);
\fill (14,-15) rectangle (15,-16);
\fill (15,-15) rectangle (16,-16);
\draw[very thin] (0,0) rectangle (16,-16);
\end{tikzpicture}\qquad
\begin{tikzpicture}[x=0.28cm,y=0.28cm]
\fill (3,-1) rectangle (4,-2);
\fill (6,-1) rectangle (7,-2);
\fill (8,-1) rectangle (9,-2);
\fill (9,-1) rectangle (10,-2);
\fill (10,-1) rectangle (11,-2);
\fill (12,-1) rectangle (13,-2);
\fill (13,-1) rectangle (14,-2);
\fill (15,-1) rectangle (16,-2);
\fill (3,-2) rectangle (4,-3);
\fill (4,-2) rectangle (5,-3);
\fill (5,-2) rectangle (6,-3);
\fill (6,-2) rectangle (7,-3);
\fill (11,-2) rectangle (12,-3);
\fill (12,-2) rectangle (13,-3);
\fill (13,-2) rectangle (14,-3);
\fill (14,-2) rectangle (15,-3);
\fill (3,-3) rectangle (4,-4);
\fill (4,-3) rectangle (5,-4);
\fill (5,-3) rectangle (6,-4);
\fill (6,-3) rectangle (7,-4);
\fill (7,-3) rectangle (8,-4);
\fill (8,-3) rectangle (9,-4);
\fill (9,-3) rectangle (10,-4);
\fill (10,-3) rectangle (11,-4);
\fill (2,-4) rectangle (3,-5);
\fill (3,-4) rectangle (4,-5);
\fill (6,-4) rectangle (7,-5);
\fill (7,-4) rectangle (8,-5);
\fill (10,-4) rectangle (11,-5);
\fill (11,-4) rectangle (12,-5);
\fill (14,-4) rectangle (15,-5);
\fill (15,-4) rectangle (16,-5);
\fill (2,-5) rectangle (3,-6);
\fill (3,-5) rectangle (4,-6);
\fill (5,-5) rectangle (6,-6);
\fill (7,-5) rectangle (8,-6);
\fill (8,-5) rectangle (9,-6);
\fill (9,-5) rectangle (10,-6);
\fill (12,-5) rectangle (13,-6);
\fill (14,-5) rectangle (15,-6);
\fill (2,-6) rectangle (3,-7);
\fill (3,-6) rectangle (4,-7);
\fill (4,-6) rectangle (5,-7);
\fill (5,-6) rectangle (6,-7);
\fill (10,-6) rectangle (11,-7);
\fill (11,-6) rectangle (12,-7);
\fill (12,-6) rectangle (13,-7);
\fill (13,-6) rectangle (14,-7);
\fill (2,-7) rectangle (3,-8);
\fill (3,-7) rectangle (4,-8);
\fill (4,-7) rectangle (5,-8);
\fill (8,-7) rectangle (9,-8);
\fill (9,-7) rectangle (10,-8);
\fill (13,-7) rectangle (14,-8);
\fill (14,-7) rectangle (15,-8);
\fill (15,-7) rectangle (16,-8);
\fill (1,-8) rectangle (2,-9);
\fill (3,-8) rectangle (4,-9);
\fill (5,-8) rectangle (6,-9);
\fill (7,-8) rectangle (8,-9);
\fill (9,-8) rectangle (10,-9);
\fill (11,-8) rectangle (12,-9);
\fill (13,-8) rectangle (14,-9);
\fill (15,-8) rectangle (16,-9);
\fill (1,-9) rectangle (2,-10);
\fill (5,-9) rectangle (6,-10);
\fill (6,-9) rectangle (7,-10);
\fill (7,-9) rectangle (8,-10);
\fill (8,-9) rectangle (9,-10);
\fill (10,-9) rectangle (11,-10);
\fill (11,-9) rectangle (12,-10);
\fill (12,-9) rectangle (13,-10);
\fill (1,-10) rectangle (2,-11);
\fill (4,-10) rectangle (5,-11);
\fill (6,-10) rectangle (7,-11);
\fill (7,-10) rectangle (8,-11);
\fill (9,-10) rectangle (10,-11);
\fill (12,-10) rectangle (13,-11);
\fill (14,-10) rectangle (15,-11);
\fill (15,-10) rectangle (16,-11);
\fill (1,-11) rectangle (2,-12);
\fill (4,-11) rectangle (5,-12);
\fill (6,-11) rectangle (7,-12);
\fill (8,-11) rectangle (9,-12);
\fill (10,-11) rectangle (11,-12);
\fill (11,-11) rectangle (12,-12);
\fill (13,-11) rectangle (14,-12);
\fill (15,-11) rectangle (16,-12);
\fill (1,-12) rectangle (2,-13);
\fill (2,-12) rectangle (3,-13);
\fill (5,-12) rectangle (6,-13);
\fill (6,-12) rectangle (7,-13);
\fill (9,-12) rectangle (10,-13);
\fill (10,-12) rectangle (11,-13);
\fill (13,-12) rectangle (14,-13);
\fill (14,-12) rectangle (15,-13);
\fill (1,-13) rectangle (2,-14);
\fill (2,-13) rectangle (3,-14);
\fill (8,-13) rectangle (9,-14);
\fill (11,-13) rectangle (12,-14);
\fill (12,-13) rectangle (13,-14);
\fill (13,-13) rectangle (14,-14);
\fill (14,-13) rectangle (15,-14);
\fill (15,-13) rectangle (16,-14);
\fill (1,-14) rectangle (2,-15);
\fill (2,-14) rectangle (3,-15);
\fill (4,-14) rectangle (5,-15);
\fill (7,-14) rectangle (8,-15);
\fill (9,-14) rectangle (10,-15);
\fill (10,-14) rectangle (11,-15);
\fill (12,-14) rectangle (13,-15);
\fill (15,-14) rectangle (16,-15);
\fill (1,-15) rectangle (2,-16);
\fill (2,-15) rectangle (3,-16);
\fill (4,-15) rectangle (5,-16);
\fill (5,-15) rectangle (6,-16);
\fill (7,-15) rectangle (8,-16);
\fill (8,-15) rectangle (9,-16);
\fill (11,-15) rectangle (12,-16);
\fill (14,-15) rectangle (15,-16);
\draw[very thin] (0,0) rectangle (16,-16);
\end{tikzpicture}
\caption{Two perfect $16\times16$ binary arrays (black $=-1$, white $=+1$; $x_1$ runs downwards,
$x_2$ to the right). Left: $(-1)^{\lfloor x_1x_2/8\rfloor+x_1\lfloor x_2/4\rfloor}$. Right: the
array of Construction~\ref{tower} at $K=2$, $n=4$.}\label{fig16}
\end{figure}
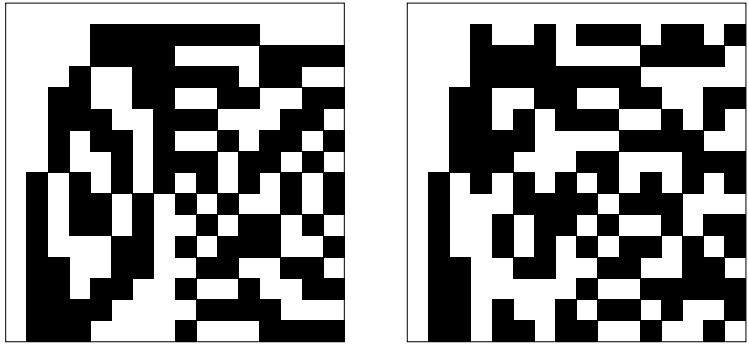

\begin{example}[$K=2$, side $16$]\label{ex16}
The array $(-1)^{\lfloor x_1x_2/8\rfloor+x_1\lfloor x_2/4\rfloor}$, $0\le x_1,x_2<16$, is a perfect
binary array. Its four $8\times8$ polyphase components have disjoint flat spectra on the cosets
of $2\Z_8^2$ (the $4\times4$ and $2\times2$ polyphase components do not), and each is a $2\times2$ twisted tiling of a $4\times4$ array of a distinct
type, i.e.\ it is an instance of Construction~\ref{thmA} with $K=2$, $M=4$. The tiles are
$$\begin{array}{cccc}
\text{type }(0,0)&\text{type }(1,0)&\text{type }(0,1)&\text{type }(1,1)\\[2pt]
\begin{pmatrix}+&+&+&+\\+&+&-&-\\+&-&+&-\\+&-&-&+\end{pmatrix}&
\begin{pmatrix}+&+&+&+\\+&+&-&-\\+&-&+&-\\+&+&-&-\end{pmatrix}&
\begin{pmatrix}+&+&-&-\\+&+&+&-\\+&-&-&+\\+&-&+&+\end{pmatrix}&
\begin{pmatrix}+&+&-&-\\+&-&+&-\\+&-&+&-\\+&+&-&-\end{pmatrix}
\end{array}$$
(a PBA, a row-wise and a column-wise quasi-perfect array, and a doubly quasi-perfect array). The
array produced by Construction~\ref{tower} at $K=2$, $n=4$ is a different perfect $16\times16$
binary array; both are shown in Figure~\ref{fig16}.
\end{example}

\begin{example}[Dillon's formula]
For $K=2$ there is an older closed formula. Dillon \cite{Dillon1990} (as stated in
\cite[Thm.~4.1]{DavisJedwab1996}) lets $r\colon\Z_{2^{a+1}}\to\{0,1\}$ send the residue $x$ to
$\lfloor x/2^a\rfloor$ and $\pi\colon\Z_{2^{a+1}}\to\Z_{2^{a+1}}$ send the residue $2^rt$, $t$ odd,
to $2^rt'$ with $tt'\equiv1\pmod{2^{a+1}}$; then $\{(x,y):r(\pi(x)y)=1\}$ is a reversible Menon
difference set in $\Z_{2^{a+1}}^2$, i.e.\ $(-1)^{r(\pi(x)y)}$ is a perfect
$2^{a+1}\times2^{a+1}$ binary array (we have checked this for $a\le5$). Its exponent is the top bit
of $\pi(x)y\bmod2^{a+1}$: a floor of a product, like the Frank--Heimiller exponent $\lfloor x_1x_2/K\rfloor$ and
Construction~V of \cite{Blake2017}, but with one factor replaced by its unit part inverted. Whether it is equivalent to the array of \eqref{lift1} at $K=2$
is open. Do Duc's \cite[Thm.~4.3]{DoDuc2021} is the $p$-ary analogue over the local ring $\Z_{p^n}$.
\end{example}

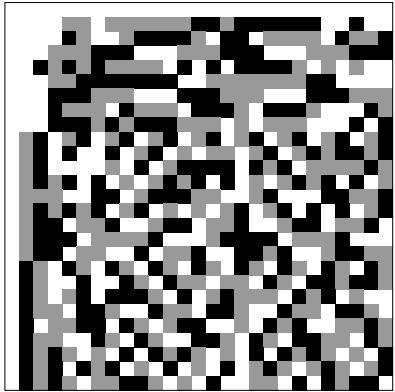
\begin{figure}[H]
\centering
\begin{tikzpicture}[x=0.19cm,y=0.19cm]
\fill[black!40] (4,-1) rectangle (5,-2);
\fill[black!40] (5,-1) rectangle (6,-2);
\fill[black!40] (7,-1) rectangle (8,-2);
\fill[black!40] (8,-1) rectangle (9,-2);
\fill[black!40] (9,-1) rectangle (10,-2);
\fill[black!40] (10,-1) rectangle (11,-2);
\fill[black!40] (11,-1) rectangle (12,-2);
\fill[black!40] (12,-1) rectangle (13,-2);
\fill[black] (13,-1) rectangle (14,-2);
\fill[black] (14,-1) rectangle (15,-2);
\fill[black!40] (15,-1) rectangle (16,-2);
\fill[black] (16,-1) rectangle (17,-2);
\fill[black] (17,-1) rectangle (18,-2);
\fill[black] (18,-1) rectangle (19,-2);
\fill[black] (19,-1) rectangle (20,-2);
\fill[black] (20,-1) rectangle (21,-2);
\fill[black] (21,-1) rectangle (22,-2);
\fill[black] (24,-1) rectangle (25,-2);
\fill[black] (4,-2) rectangle (5,-3);
\fill[black!40] (5,-2) rectangle (6,-3);
\fill[black!40] (8,-2) rectangle (9,-3);
\fill[black] (9,-2) rectangle (10,-3);
\fill[black] (10,-2) rectangle (11,-3);
\fill[black] (11,-2) rectangle (12,-3);
\fill[black] (12,-2) rectangle (13,-3);
\fill[black!40] (13,-2) rectangle (14,-3);
\fill[black] (15,-2) rectangle (16,-3);
\fill[black] (16,-2) rectangle (17,-3);
\fill[black!40] (18,-2) rectangle (19,-3);
\fill[black!40] (19,-2) rectangle (20,-3);
\fill[black!40] (20,-2) rectangle (21,-3);
\fill[black!40] (21,-2) rectangle (22,-3);
\fill[black] (23,-2) rectangle (24,-3);
\fill[black!40] (24,-2) rectangle (25,-3);
\fill[black!40] (25,-2) rectangle (26,-3);
\fill[black] (26,-2) rectangle (27,-3);
\fill[black!40] (3,-3) rectangle (4,-4);
\fill[black!40] (4,-3) rectangle (5,-4);
\fill[black!40] (5,-3) rectangle (6,-4);
\fill[black] (6,-3) rectangle (7,-4);
\fill[black] (7,-3) rectangle (8,-4);
\fill[black] (8,-3) rectangle (9,-4);
\fill[black!40] (12,-3) rectangle (13,-4);
\fill[black!40] (13,-3) rectangle (14,-4);
\fill[black!40] (14,-3) rectangle (15,-4);
\fill[black] (15,-3) rectangle (16,-4);
\fill[black] (16,-3) rectangle (17,-4);
\fill[black] (17,-3) rectangle (18,-4);
\fill[black!40] (21,-3) rectangle (22,-4);
\fill[black!40] (22,-3) rectangle (23,-4);
\fill[black!40] (23,-3) rectangle (24,-4);
\fill[black] (24,-3) rectangle (25,-4);
\fill[black] (25,-3) rectangle (26,-4);
\fill[black] (26,-3) rectangle (27,-4);
\fill[black] (2,-4) rectangle (3,-5);
\fill[black!40] (3,-4) rectangle (4,-5);
\fill[black] (4,-4) rectangle (5,-5);
\fill[black!40] (5,-4) rectangle (6,-5);
\fill[black] (6,-4) rectangle (7,-5);
\fill[black!40] (7,-4) rectangle (8,-5);
\fill[black!40] (8,-4) rectangle (9,-5);
\fill[black!40] (9,-4) rectangle (10,-5);
\fill[black!40] (10,-4) rectangle (11,-5);
\fill[black] (12,-4) rectangle (13,-5);
\fill[black] (14,-4) rectangle (15,-5);
\fill[black] (16,-4) rectangle (17,-5);
\fill[black] (17,-4) rectangle (18,-5);
\fill[black] (18,-4) rectangle (19,-5);
\fill[black] (19,-4) rectangle (20,-5);
\fill[black!40] (20,-4) rectangle (21,-5);
\fill[black!40] (22,-4) rectangle (23,-5);
\fill[black!40] (24,-4) rectangle (25,-5);
\fill[black!40] (3,-5) rectangle (4,-6);
\fill[black!40] (4,-5) rectangle (5,-6);
\fill[black] (5,-5) rectangle (6,-6);
\fill[black] (6,-5) rectangle (7,-6);
\fill[black] (7,-5) rectangle (8,-6);
\fill[black] (8,-5) rectangle (9,-6);
\fill[black] (9,-5) rectangle (10,-6);
\fill[black] (10,-5) rectangle (11,-6);
\fill[black] (11,-5) rectangle (12,-6);
\fill[black!40] (14,-5) rectangle (15,-6);
\fill[black!40] (15,-5) rectangle (16,-6);
\fill[black!40] (16,-5) rectangle (17,-6);
\fill[black!40] (17,-5) rectangle (18,-6);
\fill[black!40] (18,-5) rectangle (19,-6);
\fill[black!40] (19,-5) rectangle (20,-6);
\fill[black!40] (20,-5) rectangle (21,-6);
\fill[black] (21,-5) rectangle (22,-6);
\fill[black] (22,-5) rectangle (23,-6);
\fill[black] (3,-6) rectangle (4,-7);
\fill[black] (4,-6) rectangle (5,-7);
\fill[black] (5,-6) rectangle (6,-7);
\fill[black!40] (6,-6) rectangle (7,-7);
\fill[black!40] (7,-6) rectangle (8,-7);
\fill[black!40] (8,-6) rectangle (9,-7);
\fill[black] (12,-6) rectangle (13,-7);
\fill[black] (13,-6) rectangle (14,-7);
\fill[black] (14,-6) rectangle (15,-7);
\fill[black!40] (15,-6) rectangle (16,-7);
\fill[black!40] (16,-6) rectangle (17,-7);
\fill[black!40] (17,-6) rectangle (18,-7);
\fill[black] (21,-6) rectangle (22,-7);
\fill[black] (22,-6) rectangle (23,-7);
\fill[black] (23,-6) rectangle (24,-7);
\fill[black!40] (24,-6) rectangle (25,-7);
\fill[black!40] (25,-6) rectangle (26,-7);
\fill[black!40] (26,-6) rectangle (27,-7);
\fill[black] (3,-7) rectangle (4,-8);
\fill[black!40] (4,-7) rectangle (5,-8);
\fill[black!40] (5,-7) rectangle (6,-8);
\fill[black!40] (6,-7) rectangle (7,-8);
\fill[black] (8,-7) rectangle (9,-8);
\fill[black!40] (9,-7) rectangle (10,-8);
\fill[black!40] (10,-7) rectangle (11,-8);
\fill[black!40] (11,-7) rectangle (12,-8);
\fill[black] (13,-7) rectangle (14,-8);
\fill[black] (14,-7) rectangle (15,-8);
\fill[black] (15,-7) rectangle (16,-8);
\fill[black!40] (16,-7) rectangle (17,-8);
\fill[black] (18,-7) rectangle (19,-8);
\fill[black] (19,-7) rectangle (20,-8);
\fill[black] (20,-7) rectangle (21,-8);
\fill[black!40] (21,-7) rectangle (22,-8);
\fill[black] (25,-7) rectangle (26,-8);
\fill[black!40] (26,-7) rectangle (27,-8);
\fill[black] (3,-8) rectangle (4,-9);
\fill[black] (4,-8) rectangle (5,-9);
\fill[black] (5,-8) rectangle (6,-9);
\fill[black!40] (6,-8) rectangle (7,-9);
\fill[black] (7,-8) rectangle (8,-9);
\fill[black!40] (8,-8) rectangle (9,-9);
\fill[black] (9,-8) rectangle (10,-9);
\fill[black] (10,-8) rectangle (11,-9);
\fill[black] (11,-8) rectangle (12,-9);
\fill[black!40] (12,-8) rectangle (13,-9);
\fill[black!40] (13,-8) rectangle (14,-9);
\fill[black!40] (14,-8) rectangle (15,-9);
\fill[black!40] (16,-8) rectangle (17,-9);
\fill[black!40] (18,-8) rectangle (19,-9);
\fill[black!40] (19,-8) rectangle (20,-9);
\fill[black!40] (20,-8) rectangle (21,-9);
\fill[black] (24,-8) rectangle (25,-9);
\fill[black] (26,-8) rectangle (27,-9);
\fill[black!40] (1,-9) rectangle (2,-10);
\fill[black] (2,-9) rectangle (3,-10);
\fill[black!40] (4,-9) rectangle (5,-10);
\fill[black] (5,-9) rectangle (6,-10);
\fill[black!40] (7,-9) rectangle (8,-10);
\fill[black] (8,-9) rectangle (9,-10);
\fill[black!40] (10,-9) rectangle (11,-10);
\fill[black] (11,-9) rectangle (12,-10);
\fill[black!40] (13,-9) rectangle (14,-10);
\fill[black] (14,-9) rectangle (15,-10);
\fill[black!40] (16,-9) rectangle (17,-10);
\fill[black] (17,-9) rectangle (18,-10);
\fill[black!40] (19,-9) rectangle (20,-10);
\fill[black] (20,-9) rectangle (21,-10);
\fill[black!40] (22,-9) rectangle (23,-10);
\fill[black] (23,-9) rectangle (24,-10);
\fill[black!40] (25,-9) rectangle (26,-10);
\fill[black] (26,-9) rectangle (27,-10);
\fill[black!40] (1,-10) rectangle (2,-11);
\fill[black] (2,-10) rectangle (3,-11);
\fill[black] (4,-10) rectangle (5,-11);
\fill[black] (7,-10) rectangle (8,-11);
\fill[black!40] (9,-10) rectangle (10,-11);
\fill[black] (10,-10) rectangle (11,-11);
\fill[black!40] (12,-10) rectangle (13,-11);
\fill[black!40] (14,-10) rectangle (15,-11);
\fill[black!40] (15,-10) rectangle (16,-11);
\fill[black!40] (17,-10) rectangle (18,-11);
\fill[black] (18,-10) rectangle (19,-11);
\fill[black!40] (20,-10) rectangle (21,-11);
\fill[black] (21,-10) rectangle (22,-11);
\fill[black!40] (22,-10) rectangle (23,-11);
\fill[black] (23,-10) rectangle (24,-11);
\fill[black] (24,-10) rectangle (25,-11);
\fill[black!40] (25,-10) rectangle (26,-11);
\fill[black] (26,-10) rectangle (27,-11);
\fill[black!40] (1,-11) rectangle (2,-12);
\fill[black] (2,-11) rectangle (3,-12);
\fill[black!40] (7,-11) rectangle (8,-12);
\fill[black] (9,-11) rectangle (10,-12);
\fill[black!40] (11,-11) rectangle (12,-12);
\fill[black] (12,-11) rectangle (13,-12);
\fill[black] (13,-11) rectangle (14,-12);
\fill[black] (14,-11) rectangle (15,-12);
\fill[black] (15,-11) rectangle (16,-12);
\fill[black] (17,-11) rectangle (18,-12);
\fill[black!40] (18,-11) rectangle (19,-12);
\fill[black] (19,-11) rectangle (20,-12);
\fill[black!40] (21,-11) rectangle (22,-12);
\fill[black!40] (22,-11) rectangle (23,-12);
\fill[black!40] (23,-11) rectangle (24,-12);
\fill[black!40] (24,-11) rectangle (25,-12);
\fill[black] (25,-11) rectangle (26,-12);
\fill[black!40] (26,-11) rectangle (27,-12);
\fill[black!40] (1,-12) rectangle (2,-13);
\fill[black] (2,-12) rectangle (3,-13);
\fill[black!40] (3,-12) rectangle (4,-13);
\fill[black] (4,-12) rectangle (5,-13);
\fill[black] (6,-12) rectangle (7,-13);
\fill[black!40] (8,-12) rectangle (9,-13);
\fill[black!40] (10,-12) rectangle (11,-13);
\fill[black] (11,-12) rectangle (12,-13);
\fill[black!40] (12,-12) rectangle (13,-13);
\fill[black] (13,-12) rectangle (14,-13);
\fill[black] (15,-12) rectangle (16,-13);
\fill[black!40] (17,-12) rectangle (18,-13);
\fill[black!40] (19,-12) rectangle (20,-13);
\fill[black] (20,-12) rectangle (21,-13);
\fill[black!40] (21,-12) rectangle (22,-13);
\fill[black] (22,-12) rectangle (23,-13);
\fill[black] (24,-12) rectangle (25,-13);
\fill[black!40] (26,-12) rectangle (27,-13);
\fill[black!40] (1,-13) rectangle (2,-14);
\fill[black!40] (2,-13) rectangle (3,-14);
\fill[black!40] (3,-13) rectangle (4,-14);
\fill[black] (6,-13) rectangle (7,-14);
\fill[black] (7,-13) rectangle (8,-14);
\fill[black!40] (9,-13) rectangle (10,-14);
\fill[black] (10,-13) rectangle (11,-14);
\fill[black] (11,-13) rectangle (12,-14);
\fill[black] (12,-13) rectangle (13,-14);
\fill[black!40] (13,-13) rectangle (14,-14);
\fill[black!40] (14,-13) rectangle (15,-14);
\fill[black!40] (17,-13) rectangle (18,-14);
\fill[black] (18,-13) rectangle (19,-14);
\fill[black] (22,-13) rectangle (23,-14);
\fill[black] (23,-13) rectangle (24,-14);
\fill[black!40] (24,-13) rectangle (25,-14);
\fill[black!40] (25,-13) rectangle (26,-14);
\fill[black] (26,-13) rectangle (27,-14);
\fill[black!40] (1,-14) rectangle (2,-15);
\fill[black] (2,-14) rectangle (3,-15);
\fill[black!40] (3,-14) rectangle (4,-15);
\fill[black] (4,-14) rectangle (5,-15);
\fill[black!40] (5,-14) rectangle (6,-15);
\fill[black] (6,-14) rectangle (7,-15);
\fill[black!40] (8,-14) rectangle (9,-15);
\fill[black] (9,-14) rectangle (10,-15);
\fill[black!40] (11,-14) rectangle (12,-15);
\fill[black!40] (13,-14) rectangle (14,-15);
\fill[black!40] (15,-14) rectangle (16,-15);
\fill[black] (16,-14) rectangle (17,-15);
\fill[black!40] (18,-14) rectangle (19,-15);
\fill[black] (19,-14) rectangle (20,-15);
\fill[black] (21,-14) rectangle (22,-15);
\fill[black] (23,-14) rectangle (24,-15);
\fill[black!40] (25,-14) rectangle (26,-15);
\fill[black] (26,-14) rectangle (27,-15);
\fill[black!40] (1,-15) rectangle (2,-16);
\fill[black] (2,-15) rectangle (3,-16);
\fill[black] (3,-15) rectangle (4,-16);
\fill[black!40] (5,-15) rectangle (6,-16);
\fill[black!40] (6,-15) rectangle (7,-16);
\fill[black] (7,-15) rectangle (8,-16);
\fill[black!40] (10,-15) rectangle (11,-16);
\fill[black] (11,-15) rectangle (12,-16);
\fill[black] (12,-15) rectangle (13,-16);
\fill[black!40] (14,-15) rectangle (15,-16);
\fill[black!40] (15,-15) rectangle (16,-16);
\fill[black] (16,-15) rectangle (17,-16);
\fill[black!40] (19,-15) rectangle (20,-16);
\fill[black] (20,-15) rectangle (21,-16);
\fill[black] (21,-15) rectangle (22,-16);
\fill[black!40] (23,-15) rectangle (24,-16);
\fill[black!40] (24,-15) rectangle (25,-16);
\fill[black] (25,-15) rectangle (26,-16);
\fill[black!40] (1,-16) rectangle (2,-17);
\fill[black] (2,-16) rectangle (3,-17);
\fill[black] (3,-16) rectangle (4,-17);
\fill[black] (4,-16) rectangle (5,-17);
\fill[black!40] (6,-16) rectangle (7,-17);
\fill[black!40] (7,-16) rectangle (8,-17);
\fill[black!40] (8,-16) rectangle (9,-17);
\fill[black!40] (9,-16) rectangle (10,-17);
\fill[black] (10,-16) rectangle (11,-17);
\fill[black!40] (14,-16) rectangle (15,-17);
\fill[black] (15,-16) rectangle (16,-17);
\fill[black] (16,-16) rectangle (17,-17);
\fill[black] (17,-16) rectangle (18,-17);
\fill[black] (18,-16) rectangle (19,-17);
\fill[black!40] (20,-16) rectangle (21,-17);
\fill[black!40] (21,-16) rectangle (22,-17);
\fill[black!40] (22,-16) rectangle (23,-17);
\fill[black] (23,-16) rectangle (24,-17);
\fill[black!40] (1,-17) rectangle (2,-18);
\fill[black] (2,-17) rectangle (3,-18);
\fill[black] (3,-17) rectangle (4,-18);
\fill[black!40] (5,-17) rectangle (6,-18);
\fill[black!40] (6,-17) rectangle (7,-18);
\fill[black] (9,-17) rectangle (10,-18);
\fill[black!40] (11,-17) rectangle (12,-18);
\fill[black!40] (12,-17) rectangle (13,-18);
\fill[black] (13,-17) rectangle (14,-18);
\fill[black] (16,-17) rectangle (17,-18);
\fill[black] (17,-17) rectangle (18,-18);
\fill[black!40] (18,-17) rectangle (19,-18);
\fill[black] (19,-17) rectangle (20,-18);
\fill[black!40] (22,-17) rectangle (23,-18);
\fill[black] (23,-17) rectangle (24,-18);
\fill[black] (24,-17) rectangle (25,-18);
\fill[black!40] (25,-17) rectangle (26,-18);
\fill[black!40] (26,-17) rectangle (27,-18);
\fill[black] (1,-18) rectangle (2,-19);
\fill[black!40] (2,-18) rectangle (3,-19);
\fill[black] (4,-18) rectangle (5,-19);
\fill[black!40] (5,-18) rectangle (6,-19);
\fill[black] (7,-18) rectangle (8,-19);
\fill[black!40] (8,-18) rectangle (9,-19);
\fill[black] (10,-18) rectangle (11,-19);
\fill[black!40] (11,-18) rectangle (12,-19);
\fill[black] (13,-18) rectangle (14,-19);
\fill[black!40] (14,-18) rectangle (15,-19);
\fill[black] (16,-18) rectangle (17,-19);
\fill[black!40] (17,-18) rectangle (18,-19);
\fill[black] (19,-18) rectangle (20,-19);
\fill[black!40] (20,-18) rectangle (21,-19);
\fill[black] (22,-18) rectangle (23,-19);
\fill[black!40] (23,-18) rectangle (24,-19);
\fill[black] (25,-18) rectangle (26,-19);
\fill[black!40] (26,-18) rectangle (27,-19);
\fill[black] (1,-19) rectangle (2,-20);
\fill[black!40] (2,-19) rectangle (3,-20);
\fill[black] (5,-19) rectangle (6,-20);
\fill[black] (8,-19) rectangle (9,-20);
\fill[black!40] (9,-19) rectangle (10,-20);
\fill[black] (11,-19) rectangle (12,-20);
\fill[black!40] (12,-19) rectangle (13,-20);
\fill[black!40] (13,-19) rectangle (14,-20);
\fill[black!40] (15,-19) rectangle (16,-20);
\fill[black!40] (16,-19) rectangle (17,-20);
\fill[black] (18,-19) rectangle (19,-20);
\fill[black!40] (19,-19) rectangle (20,-20);
\fill[black] (21,-19) rectangle (22,-20);
\fill[black] (22,-19) rectangle (23,-20);
\fill[black!40] (23,-19) rectangle (24,-20);
\fill[black] (24,-19) rectangle (25,-20);
\fill[black] (25,-19) rectangle (26,-20);
\fill[black!40] (26,-19) rectangle (27,-20);
\fill[black] (1,-20) rectangle (2,-21);
\fill[black!40] (2,-20) rectangle (3,-21);
\fill[black!40] (4,-20) rectangle (5,-21);
\fill[black] (5,-20) rectangle (6,-21);
\fill[black] (7,-20) rectangle (8,-21);
\fill[black] (8,-20) rectangle (9,-21);
\fill[black] (9,-20) rectangle (10,-21);
\fill[black!40] (10,-20) rectangle (11,-21);
\fill[black] (12,-20) rectangle (13,-21);
\fill[black!40] (14,-20) rectangle (15,-21);
\fill[black] (15,-20) rectangle (16,-21);
\fill[black!40] (16,-20) rectangle (17,-21);
\fill[black!40] (17,-20) rectangle (18,-21);
\fill[black!40] (18,-20) rectangle (19,-21);
\fill[black] (20,-20) rectangle (21,-21);
\fill[black!40] (21,-20) rectangle (22,-21);
\fill[black] (22,-20) rectangle (23,-21);
\fill[black!40] (24,-20) rectangle (25,-21);
\fill[black] (1,-21) rectangle (2,-22);
\fill[black!40] (2,-21) rectangle (3,-22);
\fill[black!40] (3,-21) rectangle (4,-22);
\fill[black] (5,-21) rectangle (6,-22);
\fill[black] (6,-21) rectangle (7,-22);
\fill[black!40] (7,-21) rectangle (8,-22);
\fill[black] (10,-21) rectangle (11,-22);
\fill[black!40] (11,-21) rectangle (12,-22);
\fill[black!40] (12,-21) rectangle (13,-22);
\fill[black] (14,-21) rectangle (15,-22);
\fill[black] (15,-21) rectangle (16,-22);
\fill[black!40] (16,-21) rectangle (17,-22);
\fill[black] (19,-21) rectangle (20,-22);
\fill[black!40] (20,-21) rectangle (21,-22);
\fill[black!40] (21,-21) rectangle (22,-22);
\fill[black] (23,-21) rectangle (24,-22);
\fill[black] (24,-21) rectangle (25,-22);
\fill[black!40] (25,-21) rectangle (26,-22);
\fill[black] (1,-22) rectangle (2,-23);
\fill[black!40] (3,-22) rectangle (4,-23);
\fill[black!40] (4,-22) rectangle (5,-23);
\fill[black] (5,-22) rectangle (6,-23);
\fill[black] (6,-22) rectangle (7,-23);
\fill[black] (8,-22) rectangle (9,-23);
\fill[black!40] (9,-22) rectangle (10,-23);
\fill[black!40] (11,-22) rectangle (12,-23);
\fill[black] (12,-22) rectangle (13,-23);
\fill[black] (13,-22) rectangle (14,-23);
\fill[black!40] (16,-22) rectangle (17,-23);
\fill[black] (18,-22) rectangle (19,-23);
\fill[black!40] (19,-22) rectangle (20,-23);
\fill[black] (20,-22) rectangle (21,-23);
\fill[black!40] (23,-22) rectangle (24,-23);
\fill[black!40] (24,-22) rectangle (25,-23);
\fill[black] (25,-22) rectangle (26,-23);
\fill[black!40] (26,-22) rectangle (27,-23);
\fill[black] (1,-23) rectangle (2,-24);
\fill[black!40] (2,-23) rectangle (3,-24);
\fill[black!40] (3,-23) rectangle (4,-24);
\fill[black] (6,-23) rectangle (7,-24);
\fill[black!40] (7,-23) rectangle (8,-24);
\fill[black] (9,-23) rectangle (10,-24);
\fill[black!40] (10,-23) rectangle (11,-24);
\fill[black] (13,-23) rectangle (14,-24);
\fill[black] (14,-23) rectangle (15,-24);
\fill[black!40] (15,-23) rectangle (16,-24);
\fill[black] (17,-23) rectangle (18,-24);
\fill[black!40] (18,-23) rectangle (19,-24);
\fill[black] (20,-23) rectangle (21,-24);
\fill[black] (21,-23) rectangle (22,-24);
\fill[black!40] (22,-23) rectangle (23,-24);
\fill[black!40] (23,-23) rectangle (24,-24);
\fill[black] (25,-23) rectangle (26,-24);
\fill[black!40] (26,-23) rectangle (27,-24);
\fill[black] (1,-24) rectangle (2,-25);
\fill[black!40] (2,-24) rectangle (3,-25);
\fill[black] (3,-24) rectangle (4,-25);
\fill[black!40] (4,-24) rectangle (5,-25);
\fill[black!40] (6,-24) rectangle (7,-25);
\fill[black] (8,-24) rectangle (9,-25);
\fill[black] (10,-24) rectangle (11,-25);
\fill[black!40] (11,-24) rectangle (12,-25);
\fill[black] (12,-24) rectangle (13,-25);
\fill[black!40] (13,-24) rectangle (14,-25);
\fill[black!40] (15,-24) rectangle (16,-25);
\fill[black] (17,-24) rectangle (18,-25);
\fill[black] (19,-24) rectangle (20,-25);
\fill[black!40] (20,-24) rectangle (21,-25);
\fill[black] (21,-24) rectangle (22,-25);
\fill[black!40] (22,-24) rectangle (23,-25);
\fill[black!40] (24,-24) rectangle (25,-25);
\fill[black] (26,-24) rectangle (27,-25);
\fill[black] (1,-25) rectangle (2,-26);
\fill[black!40] (2,-25) rectangle (3,-26);
\fill[black] (3,-25) rectangle (4,-26);
\fill[black] (5,-25) rectangle (6,-26);
\fill[black!40] (6,-25) rectangle (7,-26);
\fill[black] (7,-25) rectangle (8,-26);
\fill[black!40] (9,-25) rectangle (10,-26);
\fill[black] (11,-25) rectangle (12,-26);
\fill[black!40] (13,-25) rectangle (14,-26);
\fill[black] (15,-25) rectangle (16,-26);
\fill[black!40] (17,-25) rectangle (18,-26);
\fill[black] (18,-25) rectangle (19,-26);
\fill[black!40] (19,-25) rectangle (20,-26);
\fill[black!40] (21,-25) rectangle (22,-26);
\fill[black] (22,-25) rectangle (23,-26);
\fill[black!40] (23,-25) rectangle (24,-26);
\fill[black!40] (25,-25) rectangle (26,-26);
\fill[black] (26,-25) rectangle (27,-26);
\fill[black] (1,-26) rectangle (2,-27);
\fill[black!40] (2,-26) rectangle (3,-27);
\fill[black] (3,-26) rectangle (4,-27);
\fill[black!40] (4,-26) rectangle (5,-27);
\fill[black!40] (6,-26) rectangle (7,-27);
\fill[black!40] (7,-26) rectangle (8,-27);
\fill[black] (8,-26) rectangle (9,-27);
\fill[black] (9,-26) rectangle (10,-27);
\fill[black!40] (10,-26) rectangle (11,-27);
\fill[black!40] (12,-26) rectangle (13,-27);
\fill[black] (14,-26) rectangle (15,-27);
\fill[black!40] (17,-26) rectangle (18,-27);
\fill[black!40] (18,-26) rectangle (19,-27);
\fill[black] (20,-26) rectangle (21,-27);
\fill[black] (22,-26) rectangle (23,-27);
\fill[black!40] (23,-26) rectangle (24,-27);
\fill[black] (24,-26) rectangle (25,-27);
\fill[black] (25,-26) rectangle (26,-27);
\draw[very thin] (0,0) rectangle (27,-27);
\end{tikzpicture}
\caption{The perfect $27\times27$ array over the cube roots of unity of Construction~\ref{tower} at
$K=3$, $n=3$, shown by exponent: white $=0$, grey $=1$, black $=2$ ($x_1$ runs downwards, $x_2$ to
the right).}\label{fig27}
\end{figure}

\begin{example}[$K=3$, side $27$; $K=5$, side $125$]
For odd $K$ the largest single-floor array $\omega^{\lfloor xy/q\rfloor}$ we know of is the
Frank--Heimiller array of side $K^2$; Construction~\ref{tower} at $n=3$ gives perfect arrays of side $27$ over the cube
roots of unity (Figure~\ref{fig27}) and of side $125$ over the fifth roots of unity, and $n=4,5$ give
sides $81$ and $243$ for $K=3$.
\end{example}

\begin{figure}[H]
\centering
\includegraphics[width=0.8\textwidth]{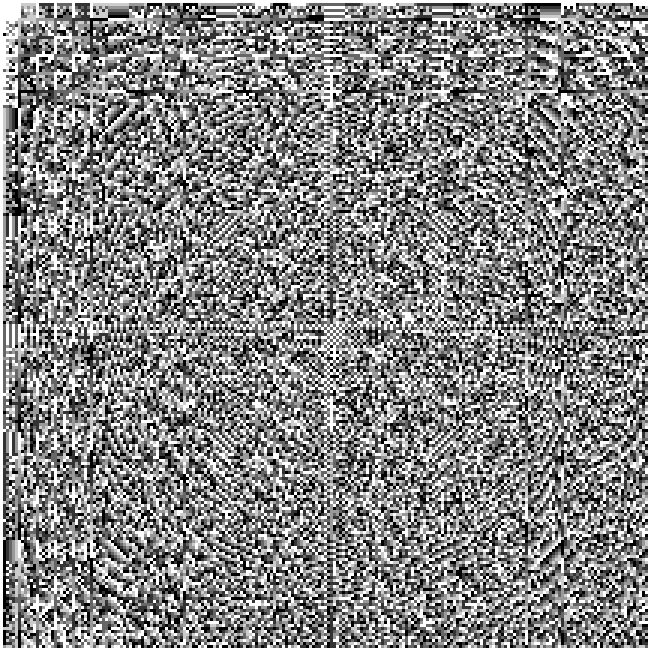}
\caption{The perfect $6^3\times6^3$ array over the sixth roots of unity of Construction~\ref{tower}
at $K=6$, $n=3$, shown by exponent on a six-level grey scale (white $=0$, black $=5$; $x_1$ runs
downwards, $x_2$ to the right).}\label{fig216}
\end{figure}

\begin{example}[$K=6$, side $216$]
Composite $K$ needs no separate treatment: Construction~\ref{tower} at $K=6$, $n=3$ gives a perfect
$216\times216$ array over the sixth roots of unity, built by the same lifts and shears as the prime
cases rather than as a Chinese-remainder product of an $8\times8$ binary and a $27\times27$ ternary
array, and using all six roots. Figure~\ref{fig216} shows it by exponent.
\end{example}

\section{Three dimensions}\label{threedsec}

Lemmas~\ref{twistlem} and \ref{interlem} are not two-dimensional: Lemma~\ref{twistlem} is stated
for an arbitrary number of coordinates with arbitrary twist orders $d_i$, and the interleaving
identity of Lemma~\ref{interlem} is a statement about one coordinate with the others carried along.
In this section we show that the lift and the shear go through in three dimensions, with one
genuinely new feature, and that for odd prime $K$ they produce perfect $K^n\times K^n\times K^n$
arrays over $\mu_K$ for all odd $n$. Types are now vectors $v\in\Z_K^3$: an array $Q$ of shape
$m=(m_1,m_2,m_3)$ over $\mu_K$ has type $v$ if $Q\rho_v$ is perfect, where
$\rho_v(x)=e^{-2\pi\sqrt{-1}\sum_i v_ix_i/(Km_i)}$. Everything below holds verbatim in $r$
dimensions; we write $r=3$ for definiteness.

\begin{theorem}[lift in three dimensions]\label{lift3thm}
Let $i\ne j$ be two of the three coordinates, let $v\in\Z_K^3$, and for each $b\in\Z_K$ let
$R^{(b)}$ be an array of shape $m$ of the type obtained from $v$ by putting $b$ in coordinate $j$.
Define $S$ of the shape obtained from $m$ by multiplying $m_i$ and $m_j$ by $K$: writing
$x_i=Ki'+a$ with $0\le a<K$ and $x_j=\bar x_j+m_jy_j$ with $0\le\bar x_j<m_j$, $0\le y_j<K$, and
leaving the third coordinate $x_k$ untouched,
$$S(x)=\omega^{\,a\,y_j}\,R^{(a)}(x\ \text{with } x_i\to i',\ x_j\to\bar x_j).$$
Then $S$ has the type obtained from $v$ by putting $0$ in coordinate $j$.
\end{theorem}
\begin{proof}
Identical to the proof of Theorem~\ref{liftthm}, with the third coordinate inert: the twisted
extension $E_{v'}[S]$ ($v'$ the claimed output type), decimated by $K$ in coordinate $i$, has
components that are twisted extensions of the $R^{(a)}$ with twist group
$\Z_K\times\Z_{K^2}\times\Z_K$ acting on coordinates $i,j,k$ and characters $(v_i,Ka,v_k)$, hence
spectra on the disjoint cosets $\{\xi_i\equiv v_i,\ \xi_j\equiv Ka\ (K^2),\ \xi_k\equiv v_k\}$,
flat by the assumed types; Lemma~\ref{interlem} in coordinate $i$ then makes $E_{v'}[S]$ coset-flat.
\end{proof}

A lift thus multiplies two of the three sides by $K$ and leaves the third alone. This is the new
feature: in two dimensions a lift multiplies both sides, so the tower stays square, whereas in
three dimensions the exponent sum $a_1+a_2+a_3$ of a shape $K^{a_1}\times K^{a_2}\times K^{a_3}$
rises by exactly $2$ at each lift, and a cube of side $K^n$ can be reached from the $K$-cube only
when $n$ is odd, through rectangular intermediate shapes. (Interleaving and twisting the same
coordinate is also allowed and multiplies that one side by $K^2$; it changes the exponent sum by
$2$ as well.)

\begin{lemma}[shear between unequal sides]\label{shear3lem}
Let $Q$ be an array of shape $m$ over $\mu_K$ of type $v$, and let $i\ne j$ be coordinates with
$m_j\mid m_i$. Then
$$Q'(x)=\omega^{\,v_j\lfloor(x_i+x_j)/m_j\rfloor}\;Q\big(x\ \text{with } x_j\to(x_i+x_j)\bmod m_j\big)$$
has type $v+(m_i/m_j)\,v_j\,e_i$. In particular it has type $v+v_je_i$ when $m_i=m_j$
(Lemma~\ref{shearlem}) and type $v$ when $K\mid m_i/m_j$.
\end{lemma}
\begin{proof}
The map $\varphi(x)=x+x_ie_j$ is an automorphism of $\prod_c\Z_{Km_c}$ precisely because
$Km_j\mid Km_i$, and $E'=E_v[Q]\circ\varphi$ satisfies $\widehat{E'}(\xi)=\widehat{E_v[Q]}(\xi-\xi_je_i)$,
which is coset-flat for the coset $v+(m_i/m_j)v_je_i+K\Z^3$. Twisted periodicity of $E_v[Q]$ gives
$E'(x+m_ie_i)=E_v[Q](\varphi(x)+m_ie_i+m_ie_j)=\omega^{v_i}\omega^{v_j(m_i/m_j)}E'(x)$, since a
shift by $m_i=(m_i/m_j)m_j$ in coordinate $j$ contributes $m_i/m_j$ twists, while shifts in the
other coordinates are unchanged; so $E'=E_{v+(m_i/m_j)v_je_i}[Q']$ with $Q'$ the restriction of $E'$
to the fundamental domain, which is the displayed formula (the floor counts the number of wraps,
now possibly more than one). Corollary~\ref{cosetflat} gives the type.
\end{proof}

So the $SL$-action of Lemma~\ref{orbitlem} is available only between coordinates of equal side,
and a shear from a side $m_j$ into a side $m_i$ with $K\mid m_i/m_j$ merely produces another array of
the same type.

\begin{lemma}[base in three dimensions]\label{base3lem}
Let $K$ be odd. The array $\omega^{x_1x_2+x_3^2}$ on $\Z_K^3$ has every axis type $(a,0,0)$ and
$(0,b,0)$, $a,b\in\Z_K$, and $\omega^{x_2x_3+x_1^2}$ has every axis type $(0,0,c)$. Consequently, for
$K$ odd, $\mathcal T_K(K,K,K)=\Z_K^3$.
\end{lemma}
\begin{proof}
With $\zeta=e^{2\pi\sqrt{-1}/K^2}$, the transform of $\omega^{x_1x_2+x_3^2}\zeta^{-ax_1}$ at
$(\alpha,\beta,\gamma)$ is
$$\sum_{x_1}\zeta^{-ax_1}\omega^{-\alpha x_1}\Big(\sum_{x_2}\omega^{x_1x_2-\beta x_2}\Big)
\Big(\sum_{x_3}\omega^{x_3^2-\gamma x_3}\Big).$$
The middle sum is $K[x_1=\beta]$, and the last is a
quadratic Gauss sum of modulus $\sqrt K$ because $K$ is odd; so the modulus is $K\sqrt K=\sqrt{K^3}$
for every $(\alpha,\beta,\gamma)$. The type $(0,b,0)$ is the same computation with the roles of
$x_1,x_2$ exchanged, and the second array is the first with coordinates permuted. All sides being
equal, the shears of Lemma~\ref{shearlem} generate $SL_3(\Z_K)$ acting on $\Z_K^3$, whose orbits are
again the classes $\gcd(v_1,v_2,v_3,K)=g$ (the proof of Lemma~\ref{orbitlem} is dimension-free), and
every class contains an axis type $(g,0,0)$.
\end{proof}

For even $K$ the term $x_3^2$ does not give a perfect sequence over $\mu_K$ (one needs $2K$-th
roots), and no perfect binary $2\times2\times2$ array exists at all, the number of entries not being
a square; the base for even $K$ is left open.

\begin{theorem}[three-dimensional tower]\label{tower3thm}
Let $K$ be an odd prime. Then $\mathcal T_K$ is all of $\Z_K^3$ at each of the shapes
$$K\times K\times K,\quad K^2\times K^2\times K,\quad K^3\times K^2\times K^2,\quad K^3\times K^3\times K^3,$$
and, repeating the three lifts, at $K^n\times K^n\times K^n$ for every odd $n$. In particular there
is a perfect $K^n\times K^n\times K^n$ array over $\mu_K$ for every odd prime $K$ and odd $n$.
\end{theorem}
\begin{proof}
The base is Lemma~\ref{base3lem}. From the full type set at $K\times K\times K$, the lifts of
Theorem~\ref{lift3thm} with $(i,j)=(1,2)$ and $(2,1)$ give, at $K^2\times K^2\times K$, every type
$(v_1,0,v_3)$ and every $(0,v_2,v_3)$; the equal-side shears between coordinates $1$ and $2$ act as
$SL_2(\Z_K)$ on $(v_1,v_2)$ with $v_3$ fixed, and since $K$ is prime every nonzero $(v_1,v_2)$ lies in
the orbit of $(1,0)$, so every type is present. From $K^2\times K^2\times K$, the lifts with
$(i,j)=(3,1)$ and $(1,3)$ give, at $K^3\times K^2\times K^2$, every $(0,v_2,v_3)$ and every
$(v_1,v_2,0)$. The shears between the equal sides $2$ and $3$ send $(v_1,v_2,0)$ to $(v_1,v_2,kv_2)$
and $(v_1,v_2,v_3)$ to $(v_1,v_2+lv_3,v_3)$, which for $K$ prime reach every $(v_1,v_2,v_3)$ with
$v_2\ne0$ or $v_3\ne0$ from the types already present, and $(v_1,0,0)$ is a lift output. From
$K^3\times K^2\times K^2$, the lifts with $(i,j)=(2,3)$ and $(3,2)$ give every $(v_1,v_2,0)$ and
$(v_1,0,v_3)$ at $K^3\times K^3\times K^3$, and now all sides are equal, so Lemma~\ref{orbitlem}
applies as in the base. The same three lifts applied to the full type set at the $K^n$-cube give
the $K^{n+2}$-cube.
\end{proof}

We have verified the theorem exactly for $K=3$ up to the $27$-cube and for $K=5$ up to the
$125$-cube, every type at every intermediate shape. The restriction to prime $K$ is only used to
fill the rectangular intermediate shapes by shears, where the full $SL_3$ action is not available.

\section{Relation to the literature}\label{litsec}

\paragraph{Binary arrays ($K=2$).}
A perfect binary array of size $s\times t$ is equivalent to a Menon--Hadamard difference set in
$\Z_s\times\Z_t$, and Jedwab's generalised perfect binary arrays of type $z\in\{0,1\}^r$
\cite[Defs.~2.7, 2.8]{Jedwab1992} are our type-$v$ arrays with $v=z$ (Definition~\ref{typedef});
his \cite[Thm.~3.2]{Jedwab1992} is the group-theoretic reading of Section~\ref{defsec}. The lift
\eqref{liftdef} at $K=2$ is \cite[Thm.~7.5]{Jedwab1992}, generalising the Jedwab--Mitchell doubling
\cite{JedwabMitchell1988,JedwabMitchell1990} (see the discussion following Theorem~\ref{liftthm}),
and the shear \eqref{shear1} at $K=2$ is \cite[Thm.~7.6]{Jedwab1992}, generalising Wild
\cite{Wild1988}. The two-dimensional versions are written out in \cite{JMPW1994}: Theorem~2.8
there (PBA + row-wise QPBA $\Rightarrow$ PBA$(2s,2t)$), Theorem~3.7 (RQPBA + DQPBA $\Rightarrow$
RQPBA$(2s,2t)$), the $c$-transform of Definition~4.4, which is the shear of the twisted extension
$[A\ A;\,-A\ -A]$ restricted to one tile, and Theorem~4.5 ($ct\equiv s\pmod{2s}$ makes the
$c$-transform carry row-wise to doubly quasi-perfect), with Wild's criterion as Corollary~4.6;
their Theorem~2.3 is the ``if and only if'' form of the one-directional interleaving criterion for
$K=2$ (interleaving is perfect iff the pieces are complementary and a cross-correlation sum
vanishes). The infinite families of quasi-perfect arrays of \cite{JedwabMitchell1990}
(Corollary~1 there: QPBA$(s,s)$ gives QPBA and DQPBA$(2^ns,2^ns)$ for all $n$), the four families of
\cite{JMPW1994} and the recursive constructions of \cite[\S9]{Jedwab1992} are built from these
operations. In particular \cite[Thm.~9.2(ii)]{Jedwab1992} gives a GPBA$(2^{a_1},\dots,2^{a_r})$ of
every nonzero type $z$ whenever $a_i\le\lceil x/2\rceil$ for $z_i=1$ and $a_i\le\lceil x/2\rceil+1$
for $z_i=0$, $x=\sum a_i$, which for $r=2$ and $a_1=a_2=n$ is our full type set at side $2^n$; and
Jedwab's recursion realises every abelian $2$-group meeting the exponent bound. The existence of PBA$(2^n,2^n)$ for all $n$ is therefore classical: necessity of
the bound $\exp(G)\le2^{d+2}$ for a Menon difference set in an abelian $2$-group of order $2^{2d+2}$
is Turyn's \cite{Turyn1965}, sufficiency is Kraemer's \cite{Kraemer1993} after Davis
\cite{Davis1991} (\cite[Thm.~2.2]{DavisJedwab1996}); the square case alone is
\cite[Thm.~4.2]{Davis1991}, obtained by transporting the $K$-matrix structure of the rectangular
group $(2^{d+2},2^d)$ \cite[Cor.~3.8]{Davis1991} to $(2^{d+1},2^{d+1})$, so Kraemer's completion is
needed only beyond rank two; and Dillon \cite{Dillon1990} gives the
explicit rank-two formula quoted above (\cite[Thm.~4.1]{DavisJedwab1996}). The Davis--Jedwab survey
also records Jedwab's result that a binary supplementary quadruple of size $s_1\times\dots\times s_r$
gives a Menon difference set in $G\times\Z_{s_1}\times\dots\times\Z_{s_r}$ for every abelian
$2$-group $G$ meeting the exponent bound (\cite[Thm.~2.7]{DavisJedwab1996}), and poses as its open
problem~3 whether an explicit construction like Dillon's exists in all abelian $2$-groups satisfying
Turyn's bound, i.e.\ beyond rank two. Nothing in
Sections~\ref{lemmasec}--\ref{towersec} is new for $K=2$ except the spectral phrasing.

\paragraph{Prime $K=p$ and relative difference sets.}
For $K=p$ our arrays are the semiregular relative difference sets of Section~\ref{defsec}, and the
groups $\Z_{p^n}^2\times\Z_p$ and $\Z_{p^{n+1}}\times\Z_{p^n}$ fall under the theorem of Ma and
Schmidt \cite{MaSchmidt1995}: a $(p^{2c},p,p^{2c},p^{2c-1})$-relative difference set exists in an
abelian group $G$ of order $p^{2c+1}$ relative to a subgroup of order $p$ if and only if
$\exp(G)\le p^{c+1}$, for every subgroup $N$ of order $p$ (the bound is Pott's \cite{Pott1994};
the constructions are Davis's \cite{Davis1992} as improved in their Theorem~2.1, the $K$-matrix
method of Davis and Kraemer in their Theorem~2.2, and the inductive Theorem~2.3, which passes from
$\langle\alpha^{p^2}\rangle\times B$ to $\langle\alpha\rangle\times B$ through one cyclic factor). So
for prime $K$ every type at every side $p^n$ was known to exist. Those constructions describe the
difference set as a union of cosets $D_{i_1,\dots,i_t}g_{i_1,\dots,i_t}$ of subgroups of a direct
factor, chosen so that each character is nonprincipal on all but one of them; what we add is the
closed index function and the uniform treatment of all $K^2$ types by one lift. The inductive
Theorem~2.3 of \cite{MaSchmidt1995} (with its siblings Theorems~4.3 and 5.4 there) is the closest
earlier relative of Theorem~\ref{liftthm}, and it is worth unwinding. Its set $R_1$ is the smaller
relative difference set $R_0\subset\langle\alpha^{p^2}\rangle\times B$ carried into
$\langle\alpha^p\rangle\times B$ by $\alpha^{p^2}\mapsto\alpha^p$ (Davis's map $P(a,b)=(2a,b)$ of
\cite[\S3]{Davis1991} is the $K=2$ version), and $\langle\alpha^{p^{e-1}}\rangle R_1$ then tiles it
$p$ times without twist; in array language this is exactly the residue-$0$ component $R^{(0)}$
of \eqref{liftdef}, tiled $K$ times untwisted along the second coordinate. The remaining $p-1$
residue classes of the cyclic coordinate are filled by Davis-type coset pieces
$D_{i_1,\dots,i_t}g_{i_1,\dots,i_t}$ with $(i_1,p)=1$ rather than by twisted tilings of typed arrays.
Since the stretch and the tiling act on the same coordinate, Theorem~2.3 is a lift by $p^2$ in one
direction whose nonzero polyphase components are supplied from outside, which is precisely the shape
that the one-dimensional discussion at the end of Section~\ref{towersec} shows a sequence tower
would need; the second factor $B$ is what makes those components available there. Related existence
results: Leung and Ma's local-ring constructions \cite{LeungMa1990}, Ma and Pott \cite{MaPott1995},
Schmidt \cite{Schmidt1997}, and Ma and Schmidt \cite{MaSchmidt2000}, who prove the unified
exponent bound (Theorem~3.2 there: for odd $p$ an abelian group of order $p^{2c+b}$ containing a
$(p^{2c},p^b,p^{2c},p^{2c-b})$-relative difference set has exponent at most $p^{c+1}$) and give a
local-ring construction: for a chain ring $R$ with maximal ideal $I$ and residue field of order $q$,
and any group $G$ of order $q^{2u+1}$ containing $I\times I$ in its centre, a
$(q^{2u},q,q^{2u},q^{2u-1})$-relative difference set in $G$ relative to $I^u\times\{0\}\subset
I\times I$ (their Theorem~4.2; Corollary~4.3 lists the groups $\Z_{p^n}^{2t}\times\Z_{p^{n-1}}^{2(s-t)}$
that can occur as $I\times I$). Hou and Sehgal \cite[Thm.~3.2]{HouSehgal2000} generalise this from
chain rings to local rings $R$ with a unique minimal ideal: with $M$ the maximal ideal, the sets
$D_a=\{(au+b(f(a)+u),u):u\in M,\ b\in B\}$ form a building set of $M\times M$ relative to
$r(M)\times\{0\}$. For the chain ring
$R=\Z_{p^{n+1}}$ this is $M\times M\cong\Z_{p^n}^2$ and $r(M)\cong\Z_p$, so the construction lands
in $\Z_{p^{n+1}}\times\Z_{p^n}$ relative to $\langle p^n\rangle\times0$ (our type $(1,0)$) but in
$\Z_{p^n}^2\times\Z_p$ only relative to a subgroup of a cyclic factor, not the direct factor (so not
our type $0$). For a Galois ring $R=\mathrm{GR}(p^{n+1},e)$ the forbidden subgroup $r(M)$ is
elementary abelian of order $p^e$, so those relative difference sets are not arrays over $p^e$-th
roots of unity; the prime-power case of Theorem~\ref{mainthm} therefore rests on the $p$-ary case
and on the inclusion $\mu_p\subset\mu_{p^e}$, as stated after Theorem~\ref{towerthm}.

\paragraph{Group-invariant Butson Hadamard matrices.}
A perfect $K$-phase array indexed by an abelian group $G$ is the same as a $G$-invariant Butson
Hadamard matrix BH$(G,K)$ \cite[Lemma~4.3]{DoDuc2020} and, for $G=\Z_q^m$, a generalised bent
function $\Z_q^m\to\Z_K$ in the sense of Kumar, Scholtz and Welch \cite{KSW1985}, and this literature settles existence for
every $K$. Do Duc's \cite[Thm.~4.3]{DoDuc2021}, applied to the local ring $R=\Z_{p^n}$, gives
BH$(\Z_{p^n}\times\Z_{p^n},p)$ for all $n$, i.e.\ perfect $p^n\times p^n$ $p$-ary arrays, by an
explicit Dillon-type formula. Schmidt, Wong and Xiang \cite[Cor.~4.3]{SWX2021} construct
BH$(\Z_{p^a}\times\Z_{p^a},h)$ whenever $p-1$ and $p$ lie in $Q(h)$, the set of non-negative
integer combinations of the prime divisors of $h$; since $p-1\notin Q(p)$ this never covers $h=p$,
and it is a source of composite-alphabet arrays on $\Z_{p^a}^2$ (for instance $6$-phase arrays for
$p=3$ or $5$) rather than of $p$-ary ones; their Corollary~7.4 (BH$(P\times\Z_p^a,h)$ for
$p\in Q(h)$) concerns a different group. The
bilinear-form construction of Do Duc and Schmidt \cite[Thm.~2.1, Cor.~2.3]{DoDucSchmidt2019}
covers BH$(\Z_{K^n}^2,K)$ for $n\le2$ (their condition is $\nu_p(K)\ge\lceil\nu_p(\exp G)/2\rceil$
for every prime $p\mid|G|$, with a parity proviso at $p=2$). The table of known two-dimensional
perfect arrays in \cite[\S4]{DoDuc2020} (Jedwab--Mitchell, Kopilovich, Wild for $h=2$; Blake, Hall
and Tirkel for $n\times n$, $n\times n^2$, $n^2\times n^2$ and $n^{2k+1}\times n^{2k+1}$ over $n$,
$2n$ or $n^{k+1}$ roots) and the new family of \cite[Thm.~4.4]{DoDuc2020}, which needs
$n_i\mid(h,n_i)^2$, all have side at most $h^2$ at alphabet $h$ except the binary families, so
none of them is a fixed-alphabet $K^n\times K^n$ family for $K>2$. Schmidt's survey
\cite{SchmidtSurvey2019} is the umbrella reference. Our
Theorem~\ref{towerthm} should be read against this background as an alternative, elementary route to
these objects that stays inside the floor-function formalism of \cite{Blake2017} and keeps track of all
$K^2$ twisted types.

\paragraph{Twisted arrays for $K>2$.}
Armario, Egan and Flannery \cite{AEF2024} define generalised perfect $K$-phase arrays of type
$z\in\{0,1\}^r$ by the expansion $\varphi'(g)=\varphi(g\bmod s)+\sum_i\lfloor g_i/s_i\rfloor$, identify
the cocycle with the product of carry cocycles $\gamma_{s_i}(j,k)=\omega^{\lfloor(j+k)/s_i\rfloor}$
\cite[Prop.~4.1]{AEF2024}, and prove for $K$ prime dividing $|G|$ that a type $z\ne0$ array is a
$(v,K,v,v/K)$-relative difference set in a factor group \cite[Thm.~4.9]{AEF2024}; see also
Armario and Flannery \cite{ArmarioFlannery2019} on quasi-orthogonal cocycles. This is
Definition~\ref{typedef} and its group-theoretic reading restricted to $v\in\{0,1\}^2$; their
Example~5.3, a $3\times3$ array over cube roots of type $(1,1)$ (we have checked that it has our
type $(1,1)$ and no other), is a representative of one of the nine types at the base level of
Lemma~\ref{baselem} for $K=3$, and their Corollary~5.5 is a Kronecker-product recursion that raises
the dimension rather than the side. We have found no
earlier use of types with entries in $\{2,\dots,K-1\}$, of the $SL_2(\Z_K)$ action on types
(Lemmas~\ref{shearlem}--\ref{orbitlem}), or of a lift of the form of Theorem~\ref{liftthm} for
$K>2$. For $K=2$ the orbit statement has a counterpart: every nonzero type is equivalent to a single
axis type, which is \cite[Cor.~7.2]{Jedwab1992} and, in difference-set language, Hughes's theorem
that a non-splitting $(4t,2,4t,2t)$-relative difference set can be assumed to have a Cartesian
forbidden subgroup \cite[Thm.~4.3]{Hughes2000}; Hughes also identifies the Jedwab cocycle
$\prod_{z_i=1}\gamma_{s_i}$, shows it is a coboundary iff no twisted side is even
\cite[Lemma~3.1]{Hughes2000}, and proves that $\tau$ is a GPBA of type $z$ iff
$\prod_{z_i=1}\gamma_{s_i}\,\partial\tau$ is a Hadamard cocycle \cite[Thm.~5.3]{Hughes2000}. The
general framework for equivalence operations on cocycles is de Launey and Flannery
\cite{deLauneyFlannery2011} and Horadam \cite{Horadam2007}.

\paragraph{Quaternary arrays ($K=4$).}
Arasu and de Launey \cite{ArasuDeLauney2001} prove that a PQA$(s,t)$ is equivalent to a
$(2st,2,2st,st)$-relative difference set in $\Z_s\times\Z_t\times\Z_4$ relative to the subgroup of
order $2$ of the $\Z_4$ factor (their Theorem~13), construct PQA$(2^n3^l,2^m3^l)$ for all $n,m\ge0$
with $-4\le n-m\le4$, $n+m\ge1$ (Theorem~6), and give inflation constructions through ``inflation
polynomials'' that multiply both sides by Mersenne and other primes (Theorems~7, 8, 18--20). Theorem~6 with
$l=0$, $n=m$ gives $4^n\times4^n$ perfect quaternary arrays for all $n$; the sizes come from
relative difference sets and Davis--Jedwab's unifying construction \cite{DavisJedwab1997} rather than from a closed
formula. Hughes \cite[Thm.~2.1]{AJC2000} makes the quaternary--binary connection explicit: $p$ on
$S$ is a PQA iff $\lambda$ on $\Z_2\times S$ is a GPBA of type $(1,0,\dots,0)$, where
$p(m)=\tfrac{1-i}{2}\big(\lambda(0,m)+i\lambda(1,m)\big)$, and deduces (Corollary~2.2) that a
PQA$(2^{b_1}3^k,2^{b_2}3^k)$ exists iff $|b_2-b_1|\le4$. In our language each $4^n\times4^n$ array
of the tower is therefore also a $2\times4^n\times4^n$ binary array of type $(1,0,0)$. Our $K=4$
tower is a different recursion from theirs (side $\times4$ per step with sixteen types), and whether
its arrays coincide with any of the Arasu--de Launey families is open.

\paragraph{Sequences and the fixed-phase question.}
In one dimension the analogous tower does not exist: for $p$ an odd prime there are no $p$-ary
perfect sequences of length $p^s$, $s\ge3$, because there is no $(p^s,p,p^s,p^{s-1})$-relative
difference set in $\Z_{p^s}\times\Z_p$ \cite[Thm.~4.6]{MaNg2009}, due to Ma and Pott
\cite{MaPott1995} (the length-$p^2$ sequence of \cite[Ex.~4.5]{MaNg2009} is the Frank sequence, the
row reading of our side-$p$ base), and for cyclic groups of prime-power
order BH$(\Z_{p^\beta},p^\alpha)$ exists iff $\beta\le2\alpha$ and $(p^\beta,p^\alpha)\ne(2,2)$
(\cite[Thm.~3.6]{DoDucSchmidt2019}, \cite[Thm.~4.8]{SchmidtSurvey2019}; cf.\
\cite{HydeKraisler2017}); the discussion at the end of Section~\ref{towersec} locates the failure
of our construction in one dimension at exactly this bound. Mow's unified construction
\cite{Mow1993,Mow1996} and his conjecture bound the length of a perfect $K$-phase sequence by
$K^2$ up to the known families. Theorem~\ref{mainthm} shows that the two-dimensional fixed-phase
problem is entirely different: the apparent ceiling $2K^2$ for single floored products (attained by Construction~V of
\cite{Blake2017}) is a property of single floored products, not of perfect arrays. Earlier engineering constructions of perfect
polyphase arrays \cite{Luke1988,BomerAntweiler1992,Antweiler1990} are tabulated in
\cite{DoDuc2020}; none gives a fixed-alphabet $K^n\times K^n$ family for general $K$.

\section{Open questions}

\begin{enumerate}
\item \emph{Non-square shapes.} Theorem~\ref{liftthm} is stated for rectangular tiles, so from a full
type set at shape $m_1\times m_2$ one obtains axis types at shape $Km_1\times Km_2$; the shears,
however, need $m_1=m_2$ (or $m_2\mid m_1$ with the appropriate carry). Lifting by $K^r$ in one
direction at a time should give shapes $K^a\times K^b$ and meet the exponent bounds of
\cite{MaSchmidt1995,MaSchmidt2000}, which for prime $K$ predict $|a-b|\le2$; in difference-set
language \cite[Thm.~2.3]{MaSchmidt1995} is such a one-direction lift with $r=2$ (see
Section~\ref{litsec}), so the question is whether its outside-supplied components can be replaced
by typed arrays. For $K=2$ the shape $2^{d+1}\times2^d$ is settled by \cite[Thm.~4.4]{MaSchmidt1995}:
a $(2^{2d+1},2,2^{2d+1},2^{2d})$-relative difference set exists in an abelian group of order
$2^{2d+2}$ iff the exponent is at most $2^{d+2}$ and the forbidden subgroup lies in a cyclic subgroup
of order $4$, so every nonzero type exists at that shape and type $0$, whose forbidden subgroup is the
direct factor $\Z_2$, does not (\cite[Lemma~3.1]{Davis1992}). The rectangular case
runs into genuinely open territory: for $p$ odd, a type-$(0,1)$ array of shape $p^{d+1}\times p^d$
over $p$-th roots is a $(p^{2d+1},p,p^{2d+1},p^{2d})$-relative difference set in
$\Z_{p^{d+1}}\times\Z_{p^{d+1}}$, and a perfect array of that shape is one in
$\Z_{p^{d+1}}\times\Z_{p^d}\times\Z_p$; by \cite[Result~1.1]{MaSchmidt2000} (Davis--Jedwab) these
are exactly the two groups left open by the exponent bound, and \cite[Thm.~5.14]{MaSchmidt2000}
shows that for $p=3$, $d=1$ the first does not exist. A rectangular lift would need a full type set
at shape $p^{d+1}\times p^d$ and therefore cannot exist in general.
\item \emph{Equivalences.} Are the arrays of Construction~\ref{tower} equivalent (under the
automorphisms of $\Z_{K^n}^2$, shifts and multiplication by characters) to Dillon's binary arrays
\cite{Dillon1990}, to Do Duc's local-ring arrays \cite{DoDuc2021}, or to the Jedwab--Mitchell
families?
\item \emph{Which perfect arrays arise?} Is every perfect $K^n\times K^n$ array over $\mu_K$ whose
$K$-polyphase components have disjoint spectral supports an interleaving of twisted extensions
(Construction~\ref{thmA})? Disjoint supports together with perfection already force each support
to be a set on which the component's transform is flat; the question is whether the supports must
be cosets of $K\Z^2$ and the components twisted tilings.
\item \emph{Higher dimensions.} Section~\ref{threedsec} reaches the cubes of odd side-exponent
for odd prime $K$. Open are: the even cubes $K^n\times K^n\times K^n$ with $n\ge4$ (the $K^2$-cube
exists as the Kronecker product of the $K^2\times K^2$ array with a Frank sequence, but at side
$K^4$ no perfect sequence over $\mu_K$ is available and the parity of the lift blocks the tower);
the base for even $K$; composite odd $K$, where the rectangular intermediate shapes may not be
filled by shears; dimensions $r\ge4$, where the same parity argument reaches the $K^n$-cube only
when $r(n-1)$ is even; and which rectangular shapes are reached in general. For $K=2$ the lift and
shear in any dimension are \cite[Thms.~7.5, 7.6]{Jedwab1992}, with the shape bounds of
\cite[Thm.~9.2]{Jedwab1992}, and the Kronecker-product recursion of \cite[Cor.~5.5]{AEF2024}
raises the dimension of type-$(1,\dots,1)$ arrays.
\end{enumerate}

\bibliographystyle{abbrv}

\begin{thebibliography}{99}

\bibitem{AJC2000} G. Hughes, ``The equivalence of certain auto-correlated quaternary and binary
arrays'', \textit{Australas. J. Combin.}, vol.~22, pp.~37--40, 2000.

\bibitem{Antweiler1990} M. Antweiler, L. B\"omer and H. D. L\"uke, ``Perfect ternary arrays'',
\textit{IEEE Trans. Inform. Theory}, vol.~36, no.~3, pp.~696--705, 1990.

\bibitem{ArasuDeLauney2001} K. T. Arasu and W. de Launey, ``Two-dimensional perfect quaternary
arrays'', \textit{IEEE Trans. Inform. Theory}, vol.~47, no.~4, pp.~1482--1493, May 2001.

\bibitem{AEF2024} J. A. Armario, R. Egan and D. L. Flannery, ``Generalized partially bent functions,
generalized perfect arrays and cocyclic Butson matrices'', \textit{Cryptogr. Commun.}, vol.~16,
pp.~323--337, 2024; arXiv:2207.05735.

\bibitem{ArmarioFlannery2019} J. A. Armario and D. L. Flannery, ``Generalized binary arrays from
quasi-orthogonal cocycles'', \textit{Des. Codes Cryptogr.}, vol.~87, pp.~2405--2417, 2019;
arXiv:1902.08812.

\bibitem{Blake2017} S. T. Blake, ``Constructions for Perfect Autocorrelation Sequences and
Multi-Dimensional Arrays'', PhD thesis, Monash University, April 2017.

\bibitem{BomerAntweiler1992} L. B\"omer and M. Antweiler, ``Perfect $N$-phase sequences and arrays'',
\textit{IEEE J. Sel. Areas Commun.}, vol.~10, no.~4, pp.~782--789, 1992.

\bibitem{Davis1991} J. A. Davis, ``Difference sets in abelian 2-groups'', \textit{J. Combin. Theory
Ser.~A}, vol.~57, pp.~262--286, 1991.

\bibitem{Davis1992} J. A. Davis, ``Construction of relative difference sets in $p$-groups'',
\textit{Discrete Math.}, vol.~103, no.~1, pp.~7--15, 1992.

\bibitem{DavisJedwab1996} J. A. Davis and J. Jedwab, ``A survey of Hadamard difference sets'', in
\textit{Groups, Difference Sets and the Monster} (K. T. Arasu et al., eds.), de Gruyter, 1996,
pp.~145--156.

\bibitem{DavisJedwab1997} J. A. Davis and J. Jedwab, ``A unifying construction for difference sets'',
\textit{J. Combin. Theory Ser.~A}, vol.~80, pp.~13--78, 1997.

\bibitem{deLauneyFlannery2011} W. de Launey and D. L. Flannery, \textit{Algebraic Design Theory},
Math. Surveys Monogr.~175, Amer. Math. Soc., 2011.

\bibitem{Dillon1990} J. F. Dillon, ``Difference sets in 2-groups'', in \textit{Finite Geometries and
Combinatorial Designs}, Contemp. Math.~111, Amer. Math. Soc., 1990, pp.~65--72.

\bibitem{DoDuc2020} T. Do Duc, ``Necessary conditions for the existence of group-invariant Butson
Hadamard matrices and a new family of perfect arrays'', \textit{Des. Codes Cryptogr.}, vol.~88,
pp.~73--90, 2020; arXiv:1903.07334.

\bibitem{DoDuc2021} T. Do Duc, ``New constructions of group-invariant Butson Hadamard matrices'',
\textit{Combinatorica}, vol.~41, pp.~625--644, 2021; arXiv:1903.09824.

\bibitem{DoDucSchmidt2019} T. Do Duc and B. Schmidt, ``Bilinear forms on finite abelian groups and
group-invariant Butson Hadamard matrices'', \textit{J. Combin. Theory Ser.~A}, vol.~166,
pp.~337--351, 2019; arXiv:1903.07310.

\bibitem{Frank1962} R. L. Frank and S. A. Zadoff, ``Phase shift pulse codes with good periodic
correlation properties'', \textit{IRE Trans. Inform. Theory}, vol.~8, no.~6, pp.~381--382, 1962.

\bibitem{Heimiller1961} R. C. Heimiller, ``Phase shift pulse codes with good periodic correlation
properties'', \textit{IRE Trans. Inform. Theory}, vol.~7, no.~4, pp.~254--257, 1961.

\bibitem{Horadam2007} K. J. Horadam, \textit{Hadamard Matrices and Their Applications}, Princeton
University Press, 2007.

\bibitem{HouSehgal2000} X.-D. Hou and S. K. Sehgal, ``Two generalized constructions of relative
difference sets'', \textit{J. Algebraic Combin.}, vol.~12, pp.~145--153, 2000.

\bibitem{Hughes2000} G. Hughes, ``Non-splitting abelian $(4t,2,4t,2t)$ relative difference sets and
Hadamard cocycles'', \textit{European J. Combin.}, vol.~21, pp.~323--331, 2000.

\bibitem{HydeKraisler2017} J. Hyde and A. Kraisler, ``Circulant $q$-Butson Hadamard matrices'',
arXiv:1701.08871, 2017.

\bibitem{Jedwab1992} J. Jedwab, ``Generalized perfect arrays and Menon difference sets'',
\textit{Des. Codes Cryptogr.}, vol.~2, pp.~19--68, 1992.

\bibitem{JedwabMitchell1988} J. Jedwab and C. J. Mitchell, ``Constructing new perfect binary arrays'',
\textit{Electron. Lett.}, vol.~24, no.~11, pp.~650--652, 1988.

\bibitem{JedwabMitchell1990} J. Jedwab and C. J. Mitchell, ``Infinite families of quasiperfect and
doubly quasiperfect binary arrays'', \textit{Electron. Lett.}, vol.~26, no.~5, pp.~294--295, 1990.

\bibitem{JMPW1994} J. Jedwab, C. J. Mitchell, F. C. Piper and P. Wild, ``Perfect binary arrays and
difference sets'', \textit{Discrete Math.}, vol.~125, pp.~241--254, 1994.

\bibitem{Kraemer1993} R. G. Kraemer, ``Proof of a conjecture on Hadamard 2-groups'', \textit{J.
Combin. Theory Ser.~A}, vol.~63, pp.~1--10, 1993.

\bibitem{KSW1985} P. V. Kumar, R. A. Scholtz and L. R. Welch, ``Generalized bent functions and their
properties'', \textit{J. Combin. Theory Ser.~A}, vol.~40, pp.~90--107, 1985.

\bibitem{LeungMa1990} K. H. Leung and S. L. Ma, ``Constructions of partial difference sets and
relative difference sets on $p$-groups'', \textit{Bull. London Math. Soc.}, vol.~22, pp.~533--539,
1990.

\bibitem{Luke1988} H. D. L\"uke, ``Sequences and arrays with perfect periodic correlation'',
\textit{IEEE Trans. Aerosp. Electron. Syst.}, vol.~24, no.~3, pp.~287--294, 1988.

\bibitem{MaNg2009} S. L. Ma and W. S. Ng, ``On non-existence of perfect and nearly perfect
sequences'', \textit{Int. J. Inform. Coding Theory}, vol.~1, no.~1, pp.~15--38, 2009.

\bibitem{MaPott1995} S. L. Ma and A. Pott, ``Relative difference sets, planar functions, and
generalized Hadamard matrices'', \textit{J. Algebra}, vol.~175, pp.~505--525, 1995.

\bibitem{MaSchmidt1995} S. L. Ma and B. Schmidt, ``On $(p^a,p,p^a,p^{a-1})$-relative difference
sets'', \textit{Des. Codes Cryptogr.}, vol.~6, pp.~57--71, 1995.

\bibitem{MaSchmidt2000} S. L. Ma and B. Schmidt, ``Relative $(p^a,p^b,p^a,p^{a-b})$-difference sets:
a unified exponent bound and a local ring construction'', \textit{Finite Fields Appl.}, vol.~6,
no.~1, pp.~1--22, 2000.

\bibitem{Mow1993} W. H. Mow, ``A Study of Correlation of Sequences'', PhD thesis, The Chinese
University of Hong Kong, 1993.

\bibitem{Mow1996} W. H. Mow, ``A new unified construction of perfect root-of-unity sequences'',
\textit{Proc. IEEE 4th ISSSTA}, vol.~3, pp.~955--959, 1996.

\bibitem{Pott1994} A. Pott, ``On the structure of abelian groups admitting divisible difference
sets'', \textit{J. Combin. Theory Ser.~A}, vol.~65, pp.~202--213, 1994.

\bibitem{Schmidt1997} B. Schmidt, ``On $(p^a,p^b,p^a,p^{a-b})$-relative difference sets'',
\textit{J. Algebraic Combin.}, vol.~6, pp.~279--297, 1997.

\bibitem{SchmidtSurvey2019} B. Schmidt, ``A survey of group invariant Butson matrices and their
relation to generalized bent functions and various other objects'', in \textit{Combinatorics and
Finite Fields} (Radon Ser. Comput. Appl. Math.~23), de Gruyter, 2019, pp.~241--251.

\bibitem{SWX2021} B. Schmidt, D. Wong and Q. Xiang, ``Constructions of Butson Hadamard matrices
invariant under abelian $p$-groups'', \textit{J. Combin. Theory Ser.~A}, vol.~181, 105433, 2021.

\bibitem{Turyn1965} R. J. Turyn, ``Character sums and difference sets'', \textit{Pacific J. Math.},
vol.~15, pp.~319--346, 1965.

\bibitem{Wild1988} P. Wild, ``Infinite families of perfect binary arrays'', \textit{Electron. Lett.},
vol.~24, no.~14, pp.~845--847, 1988.

\end{thebibliography}

\end{document}